\documentclass[12pt,leqno,twoside]{article} 
\usepackage{bezier,amsmath,amssymb,stmaryrd}
\usepackage[colorlinks, linkcolor=black, citecolor=black, urlcolor=black, hyperfootnotes=true]{hyperref}                

\input xy
\xyoption{all}
\input{rk.sty}

\renewcommand{\baselinestretch}{1.1}
\newcommand{\Niso}{$\Nl\!$-isomorphism}
\newcommand{\dDe}{\dot\De}
\newcommand{\bDes}[1]{\b\De_{#1}^\#}                 
\newcommand{\ulk}[1]{\ul{#1\!}\,}
\newcommand{\CNl}{\Cl\bby\Nl}
\DeclareMathOperator{\TTT}{{\sf T}}
\DeclareMathOperator{\LLL}{{\sf L}}
\DeclareMathOperator{\KKK}{{\sf K}}
\newcommand{\tru}{\vartriangle}
\newcommand{\trd}{\triangledown}

\DeclareMathOperator{\per}{periodic}

\begin{document}

\title{On exact functors for Heller triangulated categories}
\author{Matthias K\"unzer}
\maketitle

\bq
\begin{center}{\bf Abstract}\end{center}\vspace*{2mm}

We show certain standard constructions of the theory of Verdier triangulated categories to be valid in the Heller triangulated framework as well; viz.\ Karoubi hull, exactness of adjoints, localisation.
\eq

\renewcommand{\thefootnote}{\fnsymbol{footnote}}
\renewcommand{\thefootnote}{\arabic{footnote}}
\newdir{ >}{{}*!/-5pt/@{>}}

\begin{footnotesize}
\renewcommand{\baselinestretch}{0.7}%
\parskip0.0ex%
\tableofcontents%
\parskip1.2ex%
\renewcommand{\baselinestretch}{1.0}%
\end{footnotesize}%

\setcounter{section}{-1}

\section{Introduction}
\label{SecInt}

\subsection{Extending from Verdier to Heller}
\label{SecIntExt}

The following facts are part of the classical theory that Verdier triangulated categories.

\begin{itemize}
\item Verdier triangulated categories are stable under formation of the Karoubi hull \bfcit{BS01}.
\item The Karoubi hull construction is functorial within Verdier triangulated categories and exact functors \bfcit{BS01}.
\item Verdier triangulated categories are stable under localisation at a thick subcategory \bfcit{Ve63}.
\item Such a localisation has a universal property within Verdier triangulated categories and exact functors \bfcit{Ve63}. 
\item An adjoint functor of an exact functor is exact \bfcite{Ma83}{App.\ 2, Prop.\ 11}, \bfcite{KV87}{1.6}.
\end{itemize}

We extend these assertions somewhat to fit into the Heller triangulated setting.

\begin{itemize}
\item Heller triangulated categories are stable under formation of the Karoubi hull; cf.\ Proposition~\ref{PropKar3}.
\item The Karoubi hull construction is functorial within Heller triangulated categories and exact functors; cf.\ Proposition~\ref{PropKar4a}.
\item Closed Heller triangulated categories are stable under localisation at a thick subcategory; cf.\ Proposition~\ref{PropL5}. (Concerning closedness, see remark below.)
\item Such a localisation has a universal property within closed Heller triangulated categories and exact functors; cf.\ Proposition~\ref{PropL7}. 
\item An adjoint functor of an exact functor is exact; cf.\ Proposition~\ref{PropAd1}.
\end{itemize}

In a general Heller triangulated category, it is unknown whether there exists a cone on a given morphism. This however is true if all idempotents split \bfcite{Ku05}{Lem.\ 3.1}.
It is technically convenient to extend this assertion in the following manner. Define a Heller triangulated category to be closed if this property holds; cf.\ \bfcite{Ku07}{Def.\ 13}, Definition~\ref{DefClosed}, 
Remark~\ref{RemCl0_5}, Lemma~\ref{LemCl3}. Prove that certain constructions yield closed Heller triangulated categories or preserve closedness; cf.\ \mb{\bfcite{Ku07}{Cor.\ 21},} Proposition~\ref{PropL5}.

An exact functor between Heller triangulated categories $(\Cl,\TTT,\tht)$ and $(\Cl',\TTT',\tht')$ is a pair $(F,a)$ consisting of a subexact functor $F$ and an
isotransformation $a: \TTT F\lra F\,\TTT'\,$ such that $\tht$, $\tht'$ and $a$ are compatible; cf.\ \bfcite{KV87}{Def.\ 1.4}, Definition~\ref{DefExact}. Exactness of 
such a pair can also be characterised via $n$-triangles; cf.\ Proposition~\ref{PropC1}. The deeper reason behind that fact is that closed Heller triangulated categories can, alternatively, be defined
via sets of $n$-triangles for $n\ge 0$ with suitable properties with respect to quasicyclic and folding operations, as {\sc S.\ Thomas} informed me.

The proof of the exactness of an adjoint of an exact functor does not have to make recourse to $n$-triangles. Neither does the construction of the Heller triangulation on the Karoubi hull.
This shows the convenience of the definition of a triangulation via a tuple $\tht = (\tht_n)_{n\ge 0}$ of isomorphisms between certain shift functors, and to view the $n$-triangles as accessory, if useful; 
which is no longer the point of view taken in~\bfcit{Ku96}.

\subsection{Desirables}
\label{SecIntDes}

Still missing is a precise formulation in which sense the dual of a Heller triangulated category is again 
a Heller triangulated category, and also in which sense the constructions above are compatible with duality. Moreover, we do not treat exactness of derived functors, except implicitly, in those cases
where a derived functor can be written as a composite of an adjoint of a localisation functor, an exact functor and another localisation functor. Still missing, in the Heller triangulated
context, is furthermore the exactness of the lift of the inclusion of the heart to a functor on the bounded derived category \bfcite{BBD84}{Prop.~3.1.10}, or more generally, the 
functor $Z$ appearing in the construction of \bfcite{KV87}{Ex.\ 2.3}\,; cf.\ \bfcite{KV87}{Th.\ 3.2}.

\subsection{Notations and conventions}

We use the notations and conventions from \bfcit{Ku05}. In particular, we write composition of morphisms and functors in the natural order; viz.\ morphisms as $\lraa{f}\lraa{g} = \lraa{fg} = \lrafl{24}{f\cdot g}$
and functors as $\lraa{F}\lraa{G} = \lraa{FG} = \lrafl{24}{F\st\,G}$. Similarly for transformations. 

{\it Epic} and {\it epimorphic} are synonymous, and so are {\it monic} and {\it monomorphic.}

\section{Exact functors}

Let $(\Cl,\TTT,\tht)$, $(\Cl',\TTT',\tht')$ and $(\Cl'',\TTT'',\tht'')$ be Heller triangulated categories; cf.\ \bfcite{Ku05}{Def.\ 1.5.(ii)}.

\bq
 In \bfcite{Ku05}{Def.\ 1.5.(ii)}, we required a strictly exact functor $\Cl\lra\Cl'$ to satisfy $F\,\TTT' = \TTT F$. The adjoint functor of a strictly exact functor does not always seem to be 
 {\sf strictly} exact. Following {\sc Keller} and {\sc Vossieck}, we shall prove below that if we call a functor {\it exact,} if it satisfies the weakened condition $F\,\TTT'\iso \TTT F$ instead (and
 an accordingly modified compatibility with the Heller triangulations), then an adjoint of an exact functor is exact; cf.\ \bfcite{KV87}{1.4}. 

 Nonetheless, generally speaking, usually one deals with strictly exact functors. Hence we shall also state an extra condition of shiftcompatibility on the adjunction that ensures a shiftcompatibly 
 adjoint functor of a strictly exact functor to be strictly exact.
\eq

Given $n\ge 0$ and a transformation $G\lraa{a} G'$ between subexact additive functors $\Cl\lradoublea{G}{G'}\,\Cl'$, we denote by $G^+(\b\De_n^\#)\mrafl{25}{a^+(\b\De_n^\#)} G'^+(\b\De_n^\#)$ the 
transformation given by 
\[
\big(X (a^+(\b\De_n^\#))\big)_{\be/\al} \; :=\; X_{\be/\al} a \; :\; X_{\be/\al} G\;\lra\; X_{\be/\al} G' 
\]
for $X\in\Ob\Cl(\b\De_n^\#)$, and for $\be/\al\in\b\De_n^\#\,$, i.e.\ for $\al,\,\be\,\in\,\b\De_n$ with $\be^{-1}\le\al\le\be\le\al^{+1}$. Moreover, we denote by  
$\ulk{G^+(\b\De_n^\#)}\mrafl{35}{\ulk{a^+(\b\De_n^\#)}} \ulk{G'^+(\b\De_n^\#)}\ru{5.5}$ the induced transformation between the induced functors on the stable categories.

Sometimes, we abbreviate $(\ulk{G} \lrafl{30}{\ulk{a}} \ulk{G'}) := \big(\ulk{G^+(\b\De_n^\#)}\mrafl{35}{\ulk{a^+(\b\De_n^\#)}} \ulk{G'^+(\b\De_n^\#)}\big)$.

\begin{Definition}
\label{DefExact}\rm
\ \\
\ \\
\fbox{
\begin{tabular}{p{16cm}}
A pair $(F,a)$, consisting of an additive functor $\Cl\lraa{F}\Cl'$ and a transformation $\TTT F\lraa{a} F\,\TTT'$, is called an {\it exact pair,} or an {\it exact functor,} if the following conditions hold.

\begin{tabular}{rp{14cm}}
(1) & $a$ is an isotransformation. \\
(2) & $F$ is subexact, i.e.\ its induced functor $\h\Cl\lraa{\h F}\h\Cl'$ on the Freyd categories is exact. \\
(3) & We have
$$
(\tht_n\st \ulk{F^+(\b\De_n^\#)})\cdot\ulk{a^+(\b\De_n^\#)} \;\;=\;\; \ulk{F^+(\b\De_n^\#)}\st\tht'_n
$$
for all $n\ge 0$. \\
\end{tabular}
\end{tabular}
}

In particular, provided $\TTT F = F\,\TTT'$, then $(F,1)$ is exact if and only if $F$ is strictly exact; cf.\ \mb{\bfcite{Ku05}{Def.\ 1.5.(iii)}}. In this case, we sometimes identify $F$ and $(F,1)$.
\end{Definition}

\bq
 Calling a pair $(F,a)$ an exact functor instead of an exact pair is an abuse of notation.

 We shall not discuss whether condition~(1) is redundant; we need it for the construction of $\w Y$ in \S\ref{SecExCrit}, but that may be due to the order of our arguments.

 Condition~(3) asserts that the following cylindrical diagram commutes for all $n\ge 0$.        
 $$
 \xymatrix@C14mm{
 \ulk{\Cl^+(\b\De_n^\#)} \ar[rrr]^{\ulk{F^+(\b\De_n^\#)}}\ar@/_5mm/[ddd]_{[-]^{+1}}="phf" \ar@/^5mm/[ddd]^{[-^{+1}]}="phg" 
                    & & & \;\ulk{\Cl'^+(\b\De_n^\#)} \ar@/_5mm/[ddd]_{[-]^{+1}}="psf"\ar@/^5mm/[ddd]^{[-^{+1}]}="psg" \\
 & & & & \\
 & \ar@/^1mm/@2@<+2mm>[ur]^{1}\ar@/_1mm/@2@<-2mm>[ur]_{\ulk{a^+(\b\De_n^\#)}} & & & \\
 \ul{\Cl^+(\b\De_n^\#)}\ar[rrr]^{\ul{F^+(\b\De_n^\#)}} & & & \;\ul{\Cl'^+(\b\De_n^\#)}
 \ar@2 "phf"+<6mm,0mm>;"phg"+<-6mm,0mm>^{\tht_n\ru{-0.8}}
 \ar@2 "psf"+<6mm,0mm>;"psg"+<-6mm,0mm>^{\tht'_n\ru{-0.8}}
 }
 $$ 

 I.e., using the abbreviation just introduced, we require $X\tht_n\ulk{F}\cdot X\ulk{a} = X\ulk{F}\tht'_n$ to hold in $\ulk{\Cl'^+(\bDes{n})}$ for all $X\in\Ob\ulk{\Cl^+(\bDes{n})} = \Ob\Cl^+(\bDes{n})$. 
\eq

\begin{Definition}
\label{DefExactComp}\rm
Suppose given exact functors $(F,a)$ from $\Cl$ to $\Cl'$, and $(F',a')$ from 
$\Cl'$ to $\Cl''$. The composite of $(F,a)$ and $(F',a')$ is defined to be
\[
(F,a)\st (F',a') \= (F,a)(F',a') \; := \; (FF',\, (aF')(Fa')) \= (F\st F',\, (a\st F')\cdot (F\st a')) \; .
\]
\end{Definition}

Composition is associative.

\begin{Remark}
\label{RemD1}
If $(F,a)$ and $(F',a')$ are exact, then so is their composite $(F,a)(F',a')$.
\end{Remark}

{\it Proof.} To be able to distinguish more easily, we shall make use, from the second to the last but first step, of the notation $a\st F = aF$, $F\st F' = FF'$ etc. Given $n\ge 0$, we obtain
\[
\ba{rl}
  & \left(\tht_n\st \ulk{(F\st F')^+(\b\De_n^\#)}\right)\cdot\ulk{\big((a\st F')\cdot (F\st a')\big)^+(\b\De_n^\#)} \vspace*{2mm}\\
= & \left(\tht_n\st\ulk{F^+(\b\De_n^\#)}\st\ulk{F'^+(\b\De_n^\#)}\right)\cdot\left(\ulk{a^+(\b\De_n^\#)}\st\ulk{F'^+(\b\De_n^\#)}\right)
                                                                      \cdot\left(\ulk{F^+(\b\De_n^\#)}\st\ulk{a'^+(\b\De_n^\#)}\right) \vspace*{2mm}\\
= & \left(\ulk{F^+(\b\De_n^\#)}\st\tht'_n\st\ulk{F'^+(\b\De_n^\#)}\right)\cdot\left(\ulk{F^+(\b\De_n^\#)}\st\ulk{a'^+(\b\De_n^\#)}\right) \vspace*{2mm}\\
= & \ulk{F^+(\b\De_n^\#)}\st\ulk{F'^+(\b\De_n^\#)}\st\tht''_n \vspace*{2mm}\\
= & \ulk{(F\st F')^+(\b\De_n^\#)}\st\tht''_n \; . \vspace*{2mm}\\
\ea
\]
\qed

\begin{Definition}
\label{DefExactPeriodic}\rm
Suppose given exact functors $(F,a)$ and $(G,b)$ from $(\Cl,\TTT,\tht)$ to $(\Cl',\TTT',\tht')$. 

A transformation $F\lraa{s} G$ such that $(\TTT\st\,s)\cdot b = a\cdot (s\st\TTT')$ holds, is called {\it periodic.}
\end{Definition}

The periodicity condition requires that
\[
\xymatrix@C+2mm{
X^{+1} F\ar[r]^{X^{+1}s}\ar[d]_{Xa}^\wr & X^{+1}G\ar[d]^{Xb}_\wr \\
(XF)^{+1}\ar[r]^{(Xs)^{+1}}             & (XG)^{+1}  \\
}
\]
commute for all $X\in\Ob\Cl$.

\begin{Remark}
\label{RemExPer}\rm
Suppose given exact functors $(F,a)$, $(G,b)$ and $(H,c)$ from $\Cl$ to $\Cl'$, and periodic transformations $F\lradoublea{s}{s'} G\lraa{t} H\ru{6}\ru{-2}$.

\begin{itemize}
\item[(1)] The composite $F\lraa{s\cdot t} H$ is periodic. 
\item[(2)] The identity $F\lraa{1} F$ is periodic. 
\item[(3)] If $s$ is a periodic isotransformation from $(F,a)$ to $(G,b)$, then $s^-$ is a periodic isotransformation from $(G,b)$ to $(F,a)$. 
\item[(4)] The difference $F\lrafl{28}{s-s'} G$ of two periodic transformations is periodic. 
\item[(5)] The direct sum $(F,a)\ds (G,b) := (F\ds G,a\ds b) = (F\ds G,\smatzz{a}{0}{0}{b})$ is exact, with periodic inclusions from and periodic projections to $(F,a)$ and $(G,b)$.
\end{itemize}
\end{Remark}

\begin{Definition}
\label{DefCatExFun}\rm
Write $\,\bo\Cl,\Cl'\bc_\text{\rm ex}\,$ 
for the category of the exact functors and periodic transformations from $\,\Cl$ to $\,\Cl'\,$; cf.\ Definitions~{\rm\ref{DefExact},~\ref{DefExactPeriodic}}, Remark~{\rm\ref{RemExPer}.}

Write $\,\bo\Cl,\Cl'\bc_\text{\rm st\,ex}\,$ for the full subcategory of $\,\bo\Cl,\Cl'\bc_\text{\rm ex}\,$ of the strictly exact functors and periodic transformations from $\,\Cl$ to $\,\Cl'\,$.
\end{Definition}


\section{Idempotents and cones}
\label{SecIdCone}

Let $(\Cl,\TTT,\tht)$ be a Heller triangulated category; cf.\ \bfcite{Ku05}{Def.\ 1.5}.

\subsection{A general remark on residue classes}
\label{SecGenRes}

Concerning Frobenius categories, cf.\ e.g.\ \bfcite{Ku05}{Sec.\ A.2.3}.

\begin{Remark}
\label{RemS1}
Given a full and faithful exact functor $G : \Fl\lra\Fl\,'$ of Frobenius categories that sends all bijective objects to bijective objects.
Then the induced functor $\ulk{G} : \ulk{\Fl}\lra\ulk{\Fl}\,'$ on the classical stable categories is full and faithful.
\end{Remark}

{\it Proof.} By construction, it is full. We claim that it is faithful. Suppose given $X\lra Y$ in $\Fl$ such that 
\[
(XG\lra YG) \= (XG\lra B'\lra YG)
\]
in $\Cl'$ for some bijective object $B'$ of $\Cl'$. Choose $X\lramono B$ in $\Cl$ with $B$ bijective in $\Cl$. Since $G$ preserves pure monomorphy, $XG\lra B'$ factors over $XG\lramono BG$, whence
$XG\lra YG$ factors over $XG\lramono BG$, whence $X\lra Y$ factors over $X\lramono B$.
\qed

Suppose given weakly abelian categories $\Al$ and $\Al'$. Suppose given a subexact functor $\Al\lraa{F}\Al'\,$. Suppose given $n\ge 0$. We obtain an induced functor
\[
\ulk{F^+(\bDes{n})} \; :\; \ulk{\Al^+(\bDes{n})} \;\lra\; \ulk{\Al'^+(\bDes{n})}
\]
on the respective stable categories of $n$-pretriangles. Cf.\ \bfcite{Ku05}{\S 1.2.1.3,\,\S A.6.3}.

\begin{Remark}
\label{RemRes0}
If $F$ is full and faithful, so is $\ulk{F^+(\bDes{n})}\,$.
\end{Remark}

In particular, if $F$ is the embedding of a full subcategory, we may and will also consider $\ulk{F^+(\bDes{n})}$ to be the embedding of a full subcategory.

{\it Proof.} By \bfcite{Ku05}{Prop.\ 5.5}, both $\Al^+(\b\De_n^\#)$ and $\Al'^+(\b\De_n^\#)$ are Frobenius categories; and the full and faithful functor $F^+(\bDes{n}) : \Al^+(\b\De_n^\#)\lra \Al'^+(\b\De_n^\#)$ 
induced by $F$ preserves bijective objects, viz.\ split objects, and pure short exact sequences, viz.\ pointwise split short exact sequences. So by Remark~\ref{RemS1}, the assertion follows.
\qed

\subsection{A Heller triangulation on the Karoubi hull}
\label{SecKaroubiHeller}

Let $\h\Cl$ denote the Freyd category of $\Cl\,$; cf.\  e.g.\ \bfcite{Ku05}{§A.6.3}. We consider the full and faithful functor $\Cl\lra\h\Cl$ as an embedding of a full subcategory. 
Let $\w\Cl$ denote the full subcategory of bijectives in the abelian Frobenius category $\h\Cl$. So we have full subcategories 
\[
\Cl\;\tm\;\w\Cl\;\tm\;\h\Cl\,
\]
Since the image of $\Cl$ in $\h\Cl$ is a big enough subcategory of bijectives, 
the embedding $\Cl\hra\w\Cl$ is a Karoubi hull of $\Cl\,$; cf.\ \bfcite{Ka69}{III.II}. Cf.\ also Remark~\ref{RemAppKar1}, Lemma~\ref{LemAppKar2} -- which we will not use and argue directly instead. 

We shall give a Heller triangulation on this Karoubi hull $\w\Cl$ of $\Cl$. The Verdier triangulated version of this construction is due to {\sc Balmer} 
and {\sc Schlichting;} cf.\ \bfcite{BS01}{Th.\ 1.12}. 

As a full subcategory of bijective objects in abelian Frobenius category, the category $\w\Cl$ is weakly abelian.

The shift $\TTT$ on $\Cl$ induces a shift $\h\TTT$ on $\h\Cl$, which restricts to a shift $\w\TTT$ on $\w\Cl$.

\begin{Remark}
\label{RemKar1}
Suppose given $n\ge 0$ and $X\in\Ob\,\ulk{\w\Cl^+(\bDes{n})}\ru{-2.6}\,$. There exists $Z\in\Ob\,\ulk{\Cl^+(\bDes{n})}$ such that $X$ is isomorphic to a direct summand of $Z$ in 
$\ulk{\w\Cl^+(\bDes{n})}\,$. In other words, there \mb{exists} \mb{$Z\in\Ob\,\ulk{\Cl^+(\bDes{n})}$} and a split monomorphism $X\lraa{i} Z$ in $\ulk{\w\Cl^+(\bDes{n})}\,$.
\end{Remark}

{\it Proof.} By \bfcite{Ku05}{Prop.\ 2.6}, it suffices to prove that given \mb{$X\in\Ob\,\ulk{\w\Cl(\dDe_n)}$}, there exists \mb{$Z\in\Ob\,\ulk{\Cl(\dDe_n)}$} such that $X$ is 
isomorphic, in $\ulk{\w\Cl(\dDe_n)}$, to a direct summand of $Z$. 

It suffices to prove the existence of a split monomorphism $X\lra Z$ in $\w\Cl(\dDe_n)$ with $Z\in\Ob\Cl(\dDe_n)$.

For $i\in [1,n]$, let \mb{$Y_i\in\Ob\,\w\Cl$} be such that $X_i\ds Y_i$ is isomorphic to an object in $\Cl$. Let \mb{$Y\in\Ob\,\w\Cl(\dDe_n)$} have entry $Y_i$ at position $i$ for $1\le i\le n$ 
and the morphism from position $i$ to position $j$ be zero for $1\le i < j \le n$. The diagram $X\ds Y$ has $X$ as a summand and is isomorphic to an object in $\Cl(\dDe_n)$.
\qed

\begin{Remark}
\label{RemKar2}
Given $n\ge 0$, a diagram $X\in\Ob\,\ulk{\w\Cl^+(\bDes{n})}\,$, a split monomorphism $X\lra Z$ with $Z\,\in\,\Ob\,\ulk{\Cl^+(\bDes{n})}$ and a morphism $X\lraa{x} X'$, 
then there exists a commutative quadrangle
\[
\xymatrix{
X\ru{-1}\arspm[d]\ar[r]^x & X'\ru{-1} \arspm[d] \\
Z\ar[r]                   & Z'                  \\
}
\]
in $\ulk{\w\Cl^+(\bDes{n})}$ with $Z'\in\Ob\,\ulk{\Cl^+(\bDes{n})}\,$.

Moreover, if $X\lra X'$ is a split monomorphism, we may choose $Z\lra Z'$ to be a split monomorphism.
\end{Remark}

{\it Proof.} We form
\[
\xymatrix{
X\ru{-1}\arspm[d]_{\smatez{1}{0}}\ar[r]^x & X'\ru{-1} \arspm[d]\arspm[d]^{\smatez{1}{0}} \\
X\ds Y\ar[r]_-{\smatzz{x}{0}{0}{1}}       & X'\ds Y\; ,                                  \\
}
\]
where $X\ds Y\iso Z$. By Remark~\ref{RemKar1}, there is a split monomorphism from $X'\ds Y$ to an object $Z'$ of $\Ob\,\ulk{\Cl^+(\bDes{n})}\,$.

Moreover, if $X\lraa{x} X'$ is split monic, so is the composite $(X\ds Y\lrafl{36}{\smatzz{x}{0}{0}{1}} X'\ds Y\lra Z')$.
\qed

\begin{Construction}
\label{ConKar2_5}\rm
Given $n\ge 0$, we define $[-]^{+1}\lraa{\w\tht_n} [-^{+1}]$ on $\ulk{\w\Cl^+(\bDes{n})}$ as follows. 

Given $X\in\Ob\,\ulk{\w\Cl^+(\bDes{n})}\,$, choose a split monomorphism 
$X\lraa{i} Z$ with $Z\in\Ob\,\ulk{\Cl^+(\bDes{n})}\,$, existent by Remark~\ref{RemKar1}, and choose a retraction $p$ to $i$. Define
\[
([X]^{+1} \lrafl{30}{X\w\tht_n} [X^{+1}]) \; :=\; ([X]^{+1} \lrafl{30}{[i]^{+1}} [Z]^{+1} \lrafl{30}{Z\tht_n} [Z^{+1}] \lrafl{30}{[p^{+1}]} [X^{+1}])\; .
\]
To prove that $X\w\tht_n$ is welldefined, we shall first show that it is independent of the choice of the retraction $p$. Given $d : Z \lra X$ with $id = 0$, we have to show that $[i]^{+1} Z\tht_n [d^{+1}] = 0$.
Since $[i^{+1}]$ is monic, it suffices to show that $[i]^{+1} Z\tht_n [d^{+1}][i^{+1}] = 0$. In fact,
\[
[i]^{+1} Z\tht_n [d^{+1}][i^{+1}] \= [i]^{+1} Z\tht_n [(di)^{+1}] \= [i]^{+1} [di]^{+1} Z\tht_n \= [id]^{+1} [i^{+1}] Z\tht_n \= 0 \; ,
\]
since $di$ is in $\ulk{\Cl^+(\bDes{n})}\,$.

Now assume given another split monomorphism $X\lra Z'$ with $Z'\in\Ob\ulk{\Cl^+(\bDes{n})}\,$. By Remark~\ref{RemKar2}, we may assume that this split monomorphism factors into two
split monomorphisms $X\lraa{i} Z\lraa{i'} Z'$. Let $ip = 1$ and $i'p' = 1$. Then $(ii')(p'p) = 1$, and we may conclude
\[
\hspace*{-3mm}
[ii']^{+1}(Z'\tht_n)[(p'p)^{+1}] \= [i]^{+1}[i']^{+1}\Big(Z'\tht_n[p'^{+1}]\Big)[p^{+1}]\= [i]^{+1}[i']^{+1}\Big([p']^{+1} Z\tht_n\Big)[p^{+1}] \=  [i]^{+1} Z\tht_n[p^{+1}]\; ,
\]
since $p'$ is in $\ulk{\Cl^+(\bDes{n})}\,$. 

To show that $\tht_n$ is a transformation, we suppose given a morphism $X\lraa{f} X'$ in $\ulk{\w\Cl^+(\bDes{n})}$ and have to show that $X\w\tht_n[f^{+1}] \sollgl [f]^{+1} X'\w\tht_n$.
By Remarks~\ref{RemKar1}~and~\ref{RemKar2}, we find a commutative quadrangle
\[
\xymatrix{
X\ru{-1}\ar@{ >->}~+{|*\dir{*}}[d]_i\ar[r]^f & X'\ru{-1} \ar@{ >->}~+{|*\dir{*}}[d]^{i'} \\
Z\ar[r]^g                                    & Z'  \\
}
\]
in $\ulk{\w\Cl^+(\bDes{n})}$ with $\,Z,\, Z'\,\in\,\Ob\,\ulk{\w\Cl^+(\bDes{n})}\,$. 
Choose $p$ and $p'$ such that $ip = 1$ and $i'p' = 1$.
It suffices to show that $X\w\tht_n[f^{+1}][i'^{+1}] \sollgl [f]^{+1} X'\w\tht_n[i'^{+1}]$ by monomorphy of $[i'^{+1}]$.
In fact,
\[
\ba{clclcl}
  & X\w\tht_n[f^{+1}][i'^{+1}]                      & = & X\w\tht_n[i^{+1}][g^{+1}]                     & = & [i]^{+1}Z\tht_n[p^{+1}][i^{+1}][g^{+1}]         \\
= & [i]^{+1}Z\tht_n[(pig)^{+1}]                     & = & [i]^{+1}[pig]^{+1} Z'\tht_n                   & = & [ig]^{+1} Z'\tht_n                              \\
= & [fi']^{+1} Z'\tht_n                             & = & [fi']^{+1} [p'i']^{+1} Z'\tht_n               & = & [fi']^{+1} Z'\tht_n [(p'i')^{+1}]               \\
= & [f]^{+1} [i']^{+1} Z'\tht_n [p'^{+1}] [i'^{+1}] & = & [f]^{+1}X'\w\tht_n [i'^{+1}]\; .              &   &                                                 \\
\ea
\]
Note that $Z\w\tht_n = Z\tht_n$ for $Z\in\Ob\,\ulk{\Cl^+(\bDes{n})}\,$. 

End of construction.
\end{Construction} 

\begin{Proposition}
\label{PropKar3}
\Absit\vsp{-2}
\begin{itemize}
\item[\rm (1)] The tuple $\w\tht := (\w\tht_n)_{n\ge 0}$ is the unique Heller triangulation on $(\w\Cl,\w\TTT)$ such that the full and faithful inclusion functor $\Cl\hra\w\Cl$ is
strictly exact; \mb{cf.\ \bfcite{Ku05}{Def.\ 1.5.(i,\,ii)}}. 
\item[\rm (2)] An \mb{$n$-pretriangle} $U\in\Ob\,\ulk{\Cl^+(\bDes{n})}\ru{-2.5}$ is an $n$-triangle  with respect to $(\Cl,\TTT,\tht)$ if and only if
it is an $n$-triangle with respect to $(\w\Cl,\w\TTT,\w\tht)$.
\end{itemize}
\end{Proposition}

{\it Proof.} Ad (1). We have to show that given $m,\, n\,\ge\, 0$ and a periodic monotone map $\b\De_n\llaa{p}\b\De_m\,$, we have $\ul{p}^\#\st\w\tht_m \sollgl \w\tht_n\st\ul{p}^\#$. Let us verify this at 
$X\in\Ob\,\ulk{\w\Cl^+(\bDes{n})}\,$. Choose a split monomorphism $X\lraa{i} Z$ with $Z\in\Ob\,\ulk{\Cl^+(\bDes{n})}$ by Remark~\ref{RemKar1}. It suffices to show that
$(X\ul{p}^\#)\w\tht_m [(i\ul{p}^\#)^{+1}] \sollgl (X\w\tht_n)\ul{p}^\# [(i\ul{p}^\#)^{+1}]$. In fact, we obtain
\[
\ba{rclclcl}
(X\ul{p}^\#)\w\tht_m [(i\ul{p}^\#)^{+1}] & = & [i\ul{p}^\#]^{+1}(Z\ul{p}^\#)\tht_m & = & [i\ul{p}^\#]^{+1}(Z\tht_n)\ul{p}^\# &   & \vspace*{1mm}                                \\
                                         & = & ([i]^{+1} Z\tht_n)\ul{p}^\#         & = & (X\w\tht_n [i^{+1}])\ul{p}^\#       & = & (X\w\tht_n)\ul{p}^\# [(i\ul{p}^\#)^{+1}]\; . \\
\ea
\]

We have to show that given $n\ge 0$, we have $\ul{\ffk}_n\st\w\tht_{n+1} \sollgl \w\tht_{2n+1}\st\ul{\ffk}_n\,$. Let us verify this at \mb{$X\in\Ob\,\ulk{\w\Cl^+(\bDes{2n+1})}\,$}.  
Choose a split monomorphism $X\lraa{i} Z$ with $Z\in\Ob\,\ulk{\Cl^+(\bDes{2n+1})}$ by Remark~\ref{RemKar1}. It suffices to show that 
$(X\ul{\ffk}_n)\w\tht_{n+1} [(i\ul{\ffk}_n)^{+1}] \sollgl (X\w\tht_{2n+1})\ul{\ffk}_n [(i\ul{\ffk}_n)^{+1}]$. In fact, we obtain
\[
\ba{rclclcl}
(X\ul{\ffk}_n)\w\tht_{n+1} [(i\ul{\ffk}_n)^{+1}] & = & [i\ul{\ffk}_n]^{+1} (Z\ul{\ffk}_n)\tht_{n+1} & = &  [i\ul{\ffk}_n]^{+1} (Z\tht_{2n+1})\ul{\ffk}_n &   &                                                       \\
                                                 & = & ([i]^{+1} Z\tht_{2n+1})\ul{\ffk}_n           & = & (X\w\tht_{2n+1} [i^{+1}])\ul{\ffk}_n           & = & (X\w\tht_{2n+1})\ul{\ffk}_n [(i\ul{\ffk}_n)^{+1}]\; . \\
\ea
\]

The inclusion functor $\Cl\hra\w\Cl$ is strictly exact since it strictly commutes with shift by construction, since it is subexact because the induced functor on the Freyd categories is an equivalence, and since 
$Z\w\tht_n = Z\tht_n$ for $Z\in\Ob\,\ulk{\Cl^+(\bDes{n})}\,$.

Now suppose that both $\w\tht$ and $\w\tht'$ are Heller triangulations on $(\w\Cl,\w\TTT)$ such that $\Cl\hra\w\Cl$ is strictly exact. Suppose given $n\ge 0$ and $X\in\Ob\,\ulk{\w\Cl^+(\bDes{n})}\,$. We have to show that
$X\w\tht_n \sollgl X\w\tht'_n\,$. Choose a split monomorphism $X\lraa{i} Z$ with $Z\in\Ob\,\ulk{\Cl^+(\bDes{n})}\,$; cf.\ Remark~\ref{RemKar1}. It suffices to show that $X\w\tht_n [i^{+1}]\sollgl X\w\tht'_n [i^{+1}]$. 
In fact,
\[
X\w\tht_n [i^{+1}] \= [i]^{+1} Z\w\tht_n \= [i]^{+1} Z\tht_n \= [i]^{+1} Z\w\tht'_n \= X\w\tht'_n [i^{+1}]\; .
\]

Ad (2). Suppose given an $n$-pretriangle $U\in\Ob\,\ulk{\Cl^+(\bDes{n})}\ru{-2.5}\,$. Now $U$ is an $n$-triangle with respect to $(\Cl,\TTT,\tht)$ if and only if $U\tht_n = 1\,$, and with respect 
to $(\w\Cl,\w\TTT,\w\tht)$ if and only if $U\w\tht_n = 1\,$; cf.\ \mb{\bfcite{Ku05}{Def.\ 1.5.(ii)}}. Since $U\tht_n = U\w\tht_n\,$, these assertions are equivalent. Cf.\ also \bfcite{Ku05}{Lem.\ 3.8}.
\qed

\subsection{Functoriality of the Karoubi hull}
\label{SecUPKaroubi}

\bq
 We shall prove the universal property of the Karoubi hull directly, without making recourse to Remark~\ref{RemAppKar1} and Lemma~\ref{LemAppKar2}. 
 We will make use of the universal property and the abelianness of the Freyd category, however.
\eq

\begin{Proposition}
\label{PropKar4a}
Suppose given Heller triangulated categories $(\Cl,\TTT,\tht)$, $(\Cl',\TTT',\tht')$. Call the strictly exact inclusion functors $\KKK : \Cl \lra \w\Cl$ and $\KKK' : \Cl' \lra \w\Cl'$.

\begin{itemize}
\item[\rm (1)] Suppose given an exact functor $(F,a)$ from $\Cl$ to $\Cl'$. 

We may construct an exact functor $(\w F,\w a)$ from $\w\Cl$ to $\w\Cl'$ such that 
\[
\xymatrix@C+4mm{
\Cl\ar[r]^-{(F,a)}\ar[d]_{\KKK} & \Cl'\ar[d]^{\KKK'}  \\
\w\Cl\ar[r]^-{(\w F,\w a)}      & \w\Cl'              \\
}
\]
commutes, i.e.\ such that $(F,a)(\KKK',1) = (\KKK,1)(\w F,\w a)$, i.e.\ such that $F\st\KKK' = \KKK\st\,\w F$ and $a\st\KKK' = \KKK\st\,\w a$, i.e.\ such that $u\w F = u F$ for $u\in\Mor\Cl$
and $Z\w a = Z a$ for $Z\in\Ob\Cl$. 

The functor $\w F$ and the condition $a\st\KKK' = \KKK\st\,\w a$ uniquely determines $\w a$. If $a = 1$, then $\w a = 1$.

\item[\rm (2)] Given two exact functors $(\w F_1\,,\,\w a_1)$ and $(\w F_2\,,\,\w a_2)$ such that $(F,a)(\KKK',1) = (\KKK,1)(\w F_1\,,\,\w a_1) = (\KKK,1)(\w F_2\,,\,\w a_2)$,
there exists a unique isotransformation $\w F_1\lraisoa{\ph}\w F_2$ such that $\KKK\st\,\ph = 1$, i.e.\ such that $Z\ph = 1$ for $Z\in\Ob\Cl$. This isotransformation $\ph$ is periodic.

\item[\rm (3)] Suppose given exact functors $(F,a)$ and $(G,b)$ from $\Cl$ to $\Cl'$. Suppose given a periodic transformation $s$ from $F$ to $G$. 

Construct $(\w F,\w a)$ and $(\w G,\w b)$ as in {\rm (1)}.

There exists a unique periodic transformation $\w s$ from $\w F$ to $\w G$ such that $\KKK\st\,\w s = s\st\KKK'$, i.e.\ such that $Z\w s = Zs$ for $Z\in\Ob\Cl$. 
\end{itemize}
\end{Proposition}

{\it Proof.} Given $X\in\Ob\w\Cl$, we choose $X\lrafl{28}{i_X} Z_X \lrafl{28}{p_X} X$ in $\w\Cl$ such that $i_X\cdot p_X = 1_X$ and such that $Z_X\in\Ob\Cl$.

Moreover, choose these objects and morphisms in such a way that $Z_{X\w\TTT} = Z_X\TTT$, $i_{X\w\TTT} = i_X\w\TTT$ and $p_{X\w\TTT} = p_X\w\TTT$ for $X\in\Ob\w\Cl$.

Furthermore, if $X\in\Ob\Cl$, then choose $Z_X = X$ and $i_X = 1_X$ and $p_X = 1_X\,$.

Given $X\lraa{u} Y$ in $\w\Cl$, we let $Z_X\lraa{z_u} Z_Y$ be defined by $z_u := p_X\cdot u\cdot i_Y\,$; cf.\ Remark~\ref{RemAdd1}.

{\it Ad {\rm (1).}} Since $F$ is subexact, $\h F$ is exact. Since $W$ is a summand of 
an object in $\Cl$, also $W\h F$ is a summand of an object in $\Cl'$, hence bijective. So $\w F := \h F|_{\w\Cl}^{\w\Cl'}$ is welldefined.

We want to show that the functor $\w F$ preserves weak kernels and is therefore subexact; cf.\ Lemma~\ref{LemC0_5}. In fact, given $W\lraa{w} B\lraa{f} C$ in $\w\Cl$ such that $w$ is a weak kernel of $f$, 
we get a factorisation $w = w' i$, where $K\lramonoa{i} B$ is a kernel of $f$ in $\h\Cl$. Considering an epimorphism $P\lraepifl{34}{p} K$ in $\h\Cl$ with $P\in\Ob\w\Cl$, we obtain a factorisation $pi = p'w = p'w'i$, 
whence $p = p'w'$, whence $w'$ is epic. Since $w' \h F$ is epic and $i\h F$ is a kernel of $f\h F$, we obtain that $w\h F = w\w F$ is a weak kernel of $f\h F = f\w F$. 

The universal property of the Freyd construction yields a transformation $\h a : \h\TTT \h F \lra \h F\,\h\TTT'$. 
We let the transformation $\w a : \w\TTT \w F \lra \w F\,\w\TTT'$ be defined on $X\in\Ob\w\Cl\tm\Ob\h\Cl$ as $X\w a := X\h a$.
In particular, $Z\w a = Za$ for $Z\in\Ob\Cl$.

Given $n\ge 0$, it remains to be shown that $\ulk{\w F^+(\bDes{n})}\st\w\tht'_n \sollgl (\w\tht_n\st\ulk{\w F^+(\bDes{n})})\cdot \ulk{\w a^+(\bDes{n})}$; cf.\ Definition~\ref{DefExact}. Let us verify this at
$X\in\Ob\,\ulk{\w\Cl^+(\bDes{n})}$. Let $X\lraa{i} Z$ be a split monomorphism with \mb{$Z\in\Ob\,\ulk{\Cl^+(\bDes{n})}$}, existent by Remark~\ref{RemKar1}. It suffices to show that
\[
(X\ulk{\w F^+(\bDes{n})})\w\tht'_n \cdot [(i\ulk{\w F^+(\bDes{n})})^{+1}] \;\sollgl\; (X\w\tht_n)\ulk{\w F^+(\bDes{n})}\cdot X\ulk{\w a^+(\bDes{n})}\cdot [(i\ulk{\w F^+(\bDes{n})})^{+1}] \; . 
\]
In fact, we obtain
\[
\ba{rcl}
(X\ulk{\w F^+(\bDes{n})})\w\tht'_n \cdot [(i\ulk{\w F^+(\bDes{n})})^{+1}]
& = & [i\ulk{\w F^+(\bDes{n})}]^{+1}\cdot (Z\ulk{\w F^+(\bDes{n})})\w\tht'_n                                        \vsp{1}\\
& = & [i\ulk{\w F^+(\bDes{n})}]^{+1}\cdot (Z\ulk{F^+(\bDes{n})})\tht'_n                                             \vsp{1}\\
& \aufgl{$(F,a)$ ex.} & [i\ulk{\w F^+(\bDes{n})}]^{+1}\cdot (Z\tht_n)\ulk{F^+(\bDes{n})} \cdot Z\ulk{a^+(\bDes{n})} \vsp{1}\\
& = & [i\ulk{\w F^+(\bDes{n})}]^{+1}\cdot(Z\w\tht_n)\ulk{\w F^+(\bDes{n})}\cdot Z\ulk{\w a^+(\bDes{n})}             \vsp{1}\\
& = & [i]^{+1}\ulk{\w F^+(\bDes{n})}\cdot(Z\w\tht_n)\ulk{\w F^+(\bDes{n})}\cdot Z\ulk{\w a^+(\bDes{n})}             \vsp{1}\\
& = & (X\w\tht_n)\ulk{\w F^+(\bDes{n})}\cdot [i^{+1}]\ulk{\w F^+(\bDes{n})}\cdot Z\ulk{\w a^+(\bDes{n})}            \vsp{1}\\
& = & (X\w\tht_n)\ulk{\w F^+(\bDes{n})}\cdot X\ulk{\w a^+(\bDes{n})}\cdot [(i\ulk{\w F^+(\bDes{n})})^{+1}] \; .            \\
\ea
\]

If $a = 1$, then $\h a = 1$, so $\w a = 1$.

It remains to show that $\w a$ is uniquely determined by $\w F$ and the condition $a\st\KKK' = \KKK\st\,\w a$. In fact, given $X\in\Ob\w\Cl$, we have
\[
X\w a\cdot i_X\w F\w\TTT' \= i_X\w\TTT\w F \cdot Z_X\w a \= i_X\w\TTT\w F \cdot Z_X a \; , 
\]
and $i_X\w F\w\TTT'$ is monic.

{\it Ad {\rm (2).}} Define $\ph : \w F_1 \lraiso \w F_2$ at $X\in\Ob\w\Cl$ by
\[
\xymatrix{
Z_X F\ar[r]^{p_X\w F_1}\ar[d]_1 & X\w F_1\ar[r]^{i_X\w F_1}\ar[d]^\wr_{X\ph} & Z_X F\ar[d]^1 \\
Z_X F\ar[r]^{p_X\w F_2}         & X\w F_2\ar[r]^{i_X\w F_2}                  & Z_X F\zw{;}   \\
}
\]
cf.\ Remark~\ref{RemAdd2}. 

The tuple $\ph = (X\ph)_{X\in\Ob\w\Cl}$ is actually a transformation from $\w F_1$ to $\w F_2\,$, for given $X\lraa{u} Y$ in $\w\Cl$, we obtain
\[
\ba{clcl}
    & p_X\w F_1\cdot u\w F_1\cdot Y\ph \cdot i_Y\w F_2
& = & p_X\w F_1\cdot u\w F_1\cdot i_Y\w F_1        \vsp{1}\\
  = & p_X\w F_1\cdot(u\cdot i_Y)\w F_1             
& = & p_X\w F_1\cdot(i_X\cdot z_u)\w F_1           \vsp{1}\\
  = & p_X\w F_1\cdot i_X\w F_1\cdot z_u\w F_1      
& = & ((p_X\cdot i_X) \cdot z_u) F                 \vsp{1}\\
  = & (z_u\cdot (p_Y\cdot i_Y)) F                  
& = & z_u \w F_2\cdot p_Y\w F_2\cdot i_Y\w F_2     \vsp{1}\\
  = & (z_u\cdot p_Y)\w F_2\cdot i_Y\w F_2          
& = & (p_X\cdot u)\w F_2\cdot i_Y\w F_2            \vsp{1}\\
  = & p_X\w F_2 \cdot u\w F_2\cdot i_Y\w F_2       
& = & p_X\w F_1\cdot X\ph \cdot u\w F_2\cdot i_Y\w F_2\;, \\
\ea
\] 
and $p_X\w F_1$ is epic and $i_Y\w F_2$ is monic.

Note that commutativity of the diagram above is also necessary, for we require $\KKK\st\,\ph = 1$. This ensures uniqueness of $\ph$.

It remains to show that $\ph$ is a periodic transformation from $(\w F_1\,,\,\w a_1)$ to $(\w F_2\,,\,\w a_2)$. In fact, given $X\in\Ob\w\Cl$, we get
\[
\ba{clcl}
    & X\w a_1\cdot X\ph\w\TTT' \cdot i_X \w F_2\w\TTT'
& = & X\w a_1\cdot i_X \w F_1\w\TTT'                   \vsp{1}\\
  = & i_X \w\TTT\w F_1\cdot Z_X\w a_1                  
& = & i_{X\w\TTT}\w F_1\cdot Z_X\w a_1                 \vsp{1}\\
  = & i_{X\w\TTT}\w F_1\cdot Z_X a                     
& = & i_{X\w\TTT}\w F_1\cdot Z_X\w a_2                 \vsp{1}\\
  = & X\w\TTT\ph\cdot i_{X\w\TTT}\w F_2\cdot Z_X\w a_2 
& = & X\w\TTT\ph\cdot i_X \w\TTT\w F_2\cdot Z_X\w a_2  \vsp{1}\\
  = & X\w\TTT\ph\cdot X\w a_2 \cdot i_X \w F_2\w\TTT'\; ,     \\
\ea
\]
and $i_X \w F_2\w\TTT'$ is monic.
 
{\it Ad {\rm (3).}} Define $\w s : \w F \lra \w G$ at $X\in\Ob\w\Cl$ by
\[
\xymatrix{
Z_X F\ar[r]^{p_X\w F}\ar[d]_{Z_X s} & X\w F\ar[r]^{i_X\w F}\ar[d]^\wr_{X\w s} & Z_X F\ar[d]^{Z_X s} \\
Z_X G\ar[r]^{p_X\w G}               & X\w G\ar[r]^{i_X\w G}                   & Z_X G\zw{;}         \\
}
\]
cf.\ Remark~\ref{RemAdd2}. 

The tuple $s = (Xs)_{X\in\Ob\w\Cl}$ is actually a transformation from $\w F$ to $\w G$, for given $X\lraa{u} Y$ in $\w\Cl$, we obtain
\[
\ba{clcl}
    & p_X\w F\cdot u\w F\cdot Y\w s \cdot i_Y\w G
& = & p_X\w F\cdot u\w F\cdot i_Y\w F\cdot Z_Y s        \vsp{1}\\
  = & p_X\w F\cdot(u\cdot i_Y)\w F\cdot Z_Y s            
& = & p_X\w F\cdot(i_X\cdot z_u)\w F \cdot Z_Y s        \vsp{1}\\
  = & p_X\w F\cdot i_X\w F\cdot z_u\w F \cdot Z_Y s     
& = & ((p_X\cdot i_X) \cdot z_u) F\cdot Z_Y s           \vsp{1}\\
  = & (z_u\cdot (p_Y\cdot i_Y)) F\cdot Z_Y s                 
& = & Z_X s\cdot (z_u\cdot (p_Y\cdot i_Y)) G            \vsp{1}\\ 
  = & Z_X s\cdot z_u \w G\cdot p_Y\w G\cdot i_Y\w G     
& = & Z_X s\cdot(z_u\cdot p_Y)\w G\cdot i_Y\w G         \vsp{1}\\
  = & Z_X s\cdot(p_X\cdot u)\w G\cdot i_Y\w G            
& = & Z_X s\cdot p_X\w G \cdot u\w G\cdot i_Y\w G       \vsp{1}\\
  = & p_X\w F\cdot X\w s \cdot u\w G\cdot i_Y\w G\; ,          \\
\ea
\] 
and $p_X\w F$ is epic and $i_Y\w G$ is monic.

Note that commutativity of the diagram above is also necessary, for we require $\KKK\st\,\w s = s\st\KKK'$. This ensures uniqueness of $s$.

It remains to show that $\w s$ is a periodic transformation from $(\w F,\w a)$ to $(\w G,\w b)$. In fact, given $X\in\Ob\w\Cl$, we get
\[
\ba{clcl}
    & X\w a\cdot X\w s\w\TTT' \cdot i_X \w G\,\w\TTT'
& = & X\w a\cdot i_X \w F\,\w\TTT'\cdot Z_X s\w\TTT'                 \vsp{1}\\
  = & i_X \w\TTT\w F\cdot Z_X\w a\cdot Z_X s\w\TTT'                  
& = & i_{X\w\TTT}\w F\cdot Z_X\w a\cdot Z_X s\w\TTT'                 \vsp{1}\\
  = & i_{X\w\TTT}\w F\cdot Z_X a\cdot Z_X s\TTT'                     
& = & i_{X\w\TTT}\w F\cdot Z_X\TTT s\cdot Z_X b                      \vsp{1}\\
  = & i_{X\w\TTT}\w F\cdot Z_X\w\TTT s\cdot Z_X\w b                  
& = & i_{X\w\TTT}\w F\cdot Z_{X\w\TTT} s\cdot Z_X\w b                \vsp{1}\\
  = & X\w\TTT\w s\cdot i_{X\w\TTT}\w G\cdot Z_X\w b 
& = & X\w\TTT\w s\cdot i_X \w\TTT\w G\cdot Z_X\w b                   \vsp{1}\\
  = & X\w\TTT\w s\cdot X\w b \cdot i_X \w G\,\w\TTT'\; ,                    \\
\ea
\]
and $i_X \w G\,\w\TTT'$ is monic.
\qed

\subsection{Closed Heller triangulated categories}
\label{SecClosed}

Recall that given a Heller triangulated category $(\Cl,\TTT,\tht)$, its Karoubi hull $\w\Cl$ is Heller triangulated, too; cf.\ Proposition~\ref{PropKar3}.(1). More precisely, $(\w\Cl,\w T,\w\tht)$ is Heller 
triangulated, where $\w T$ and $\w\tht$ are as in \S\ref{SecKaroubiHeller}.

\begin{Definition}
\label{DefClosed}\rm
\ \\
\ \\
\fbox{
\begin{tabular}{p{16cm}}
A Heller triangulated category $(\Cl,\TTT,\tht)$ is called {\it closed\,} if whenever $(X,Y,\w Z)$ is a $2$-triangle in $\w\Cl$ and $X,\, Y\,\in\,\Ob\,\Cl$, then $\w Z$ is isomorphic to an object of $\Cl$.
\end{tabular}
}
\end{Definition}

Cf.\ \bfcite{Ku05}{Def.\ 1.5.(i,\,iii)}.

\bq
 I do not know an example of a non-closed Heller triangulated category.
\eq

As usual, we will call $\w Z$ the {\it cone} of $X\lra Y$, being unique up to isomorphism. Thus we may rephrase that by definition, $(\Cl,\TTT,\tht)$ is closed if it is closed under taking cones in 
the Karoubi hull $\w\Cl$.

\begin{Remark}
\label{RemCl0_5}
The Heller triangulated category $(\Cl,\TTT,\tht)$ is closed if and only if given $X\lraa{f} Y$ in $\Cl$, there exists a $2$-triangle $X\lraa{f} Y\lra Z \lra X^{+1}$ in $\Cl$.
\end{Remark}

Cf.\ \bfcite{Ku07}{Def.\ A.6}.

{\it Proof.} If $(\Cl,\TTT,\tht)$ is closed, then given $X\lraa{f} Y$ in $\Cl$, there exists a $2$-triangle $X\lraa{f} Y\lra \w Z \lra X^{+1}\ru{5}$ in $\w\Cl$ by \bfcite{Ku05}{Lem.\ 3.1}, and we may 
substitute $\w Z$ isomorphically by an object $Z$ in $\Ob\,\Cl$, so we are done by \bfcite{Ku05}{Lem.\ 3.4.(4)}. 

Conversely, if we dispose of this existence property, and if we are given a $2$-triangle $(X,Y,\w Z)$ 
in $\w\Cl$ with $X,\, Y\,\in\,\Ob\,\Cl$, then there exists a $2$-triangle $(X,Y,Z)$ with $Z\in\Ob\,\Cl$, too, and we may apply \bfcite{Ku05}{Lem.\ 3.4.(6)} to conclude that $Z\iso\w Z$. 
So $(\Cl,\TTT,\tht)$ is closed.\qed

\begin{Remark}
\label{RemCl1}
If idempotents split in $\Cl$, then $(\Cl,\TTT,\tht)$ is closed.
\end{Remark}

{\it Proof.} If idempotents split in $\Cl$, then $\Cl = \w\Cl$. \qed

\begin{Remark}
\label{RemCl2}
Suppose given Heller triangulated categories $(\Cl,\TTT,\tht)$, $(\Cl',\TTT',\tht')$ and a full and faithful strictly exact functor $\Cl\lraa{F}\Cl'$. Furthermore, suppose that whenever given a 
$2$-triangle $(XF,YF,Z')$ in $\Cl'$, where $X,\, Y\,\in\, \Ob\,\Cl$, then there exists $Z\in\Ob\,\Cl$ such that $Z'\iso ZF$.

Suppose that $\Cl'$ is closed. Then $\Cl$ is closed.
\end{Remark}

{\it Proof.} Suppose given $X\lraa{f} Y$ in $\Cl$. There exists a $2$-triangle $XF\lraa{fF} YF\lra Z'\lra XF^{+1}$ in~$\Cl'$. By assumption, there exists $Z\in\Ob\Cl$ such that $ZF\iso Z'$. By isomorphic
substitution and fullness of $F$, we obtain a $2$-triangle $XF\lraa{fF} YF\lraa{gF} ZF\lraa{hF} XF^{+1}$ in $\Cl'\,$; cf.\ \bfcite{Ku05}{Lem.~3.4.(4)}. Since
\[
(X,Y,Z)\tht_2 \ulk{F^+(\bDes{2})} \= (X,Y,Z)\ulk{F^+(\bDes{2})}\tht'_2 \= (XF,YF,ZF)\tht'_2 \= 1\; ,
\]
we conclude by faithfulness of $\ulk{F^+(\bDes{2})}$ that $(X,Y,Z)\tht_2 = 1\,$; cf.\ Remark~\ref{RemRes0}, \bfcite{Ku05}{Def.\ 1.5.(ii)}.
So we are done by Remark~\ref{RemCl0_5}.
\qed

\begin{Remark}
\label{RemCl2_5}
A closed Heller triangulated category is Verdier triangulated.
\end{Remark}

{\it Proof.} Its Karoubian hull is Verdier triangulated \bfcite{Ku05}{Prop.~3.6}. An additive shift-closed sub\-category of a Verdier triangulated category that is closed under forming cones is
Verdier triangulated.\qed

\begin{Definition}
\label{DefCl2_7}\rm
Suppose given a closed Heller triangulated category $(\Cl,\TTT,\tht)$. 

Suppose given $n\ge 0$ and $Y\in\Ob\,\Cl(\dDe_n)$ and $X\in\Ob\,\Cl^{+,\,\tht=1}(\bDes{n})$ such that $X|_{\dDe_n} = Y$.

Then $Y$ is called the {\it base} of the $n$-triangle $X$.
\end{Definition}

\begin{Lemma}
\label{LemCl3}
Suppose given a closed Heller triangulated category $(\Cl,\TTT,\tht)$ and $n\ge 0$. The restriction functor $\Cl^{+,\,\tht=1}(\bDes{n})\;\lrafl{35}{(-)|_{\dDe_n}}\;\Cl(\dDe_n)\ru{6}$ is strictly 
dense, i.e.\ surjective on objects. In other words, each object $Y\in\Ob\,\Cl(\dDe_n)$ is the base of an $n$-triangle.
\end{Lemma}

\bq
 So weakening the assumption in \bfcite{Ku05}{Lem.\ 3.1} that idempotents be split in $\Cl$ to the assumption that $\Cl$ be closed, we nonetheless obtain the conclusion of loc.\ cit.
\eq

{\it Proof.} Suppose given $Y\in\Ob\,\Cl(\dDe_n)$. By \bfcite{Ku05}{Lem.\ 3.1}, we obtain an $n$-triangle \mb{$\w X\in\Ob\,\Cl^{+,\,\tht=1}(\bDes{n})$} such that $\w X|_{\dDe_n} = Y$.

By \bfcite{Ku05}{Lem.\ 3.4.(1)}, we have a triangle $(\w X_{\al/0}\,,\,\w X_{\be/0}\,,\,\w X_{\be/\al})$ whenever $0 < \al < \be < 0^{+1}$. Since $\Cl$ is closed, $\w X_{\be/\al}$ is isomorphic to an object 
of $\Cl$. Isomorphic substitution, which is permitted without leaving $\w\Cl^{+,\,\tht=1}(\bDes{n})$ by \bfcite{Ku05}{Lem.\ 3.4.(4)}, yields an $n$-triangle in $\Cl^{+,\,\tht=1}(\bDes{n})$ that 
restricts to $Y$ on $\dDe_n\,$; cf.\ Proposition~\ref{PropKar3}.(2).\qed


\section{Heller triangulated subcategories}

\begin{Definition}
\label{DefS1_6}\rm
Given a Heller triangulated category $(\Cl',\TTT',\tht')$, a full subcategory $\Cl\tm\Cl'$ is called a full {\it Heller triangulated subcategory} of $\Cl'$ if there exist $\TTT$ and $\tht$ such
that $(\Cl,\TTT,\tht)$ is a Heller triangulated category and such that the inclusion functor $\Cl\hra\Cl'$ is strictly exact. 
\end{Definition}

We remark that in this case, the automorphism $\TTT$ and the tuple of 
transformations $\tht$ are uniquely determined by $(\Cl',\TTT',\tht')$ as respective restrictions; cf.\ \bfcite{Ku05}{Def.\ 1.5.(iii)}, Remark~\ref{RemRes0}. 

\begin{Example}
\label{ExS1_7}\rm
Let $(\Cl,\TTT,\tht)$ be a Heller triangulated category. Let $\w\Cl$ be the Karoubi hull of $\Cl$, and let $(\w\Cl,\w\TTT,\w\tht_n)$ be the Heller triangulated category from Construction~\ref{ConKar2_5}. By 
Proposition~\ref{PropKar3}.(1), $\Cl$ is a Heller triangulated subcategory of $\w\Cl$.
\end{Example}

\begin{Lemma}
\label{LemS2} 
Suppose given a closed Heller triangulated category $(\Cl',\TTT',\tht')$, and a full subcategory $\Cl\tm\Cl'$ such that the following conditions {\rm (1,\,2)} hold.
\begin{itemize}
\item[{\rm (1)}] $\Cl\TTT' = \Cl$.
\item[{\rm (2)}] Given a $2$-triangle $(X,Y,Z')$ in $\Cl'$ with $X,\, Y\,\in\, \Ob\,\Cl$, then $Z'$ is isomorphic to an object of $\Cl$.
\end{itemize}
Then $\Cl$, equipped with the shift $\TTT$ and the tuple $\tht$ obtained by restriction from $\TTT'$ and $\tht'$, respectively, is a Heller triangulated subcategory of $\Cl'$. Moreover, $(\Cl,\TTT,\tht)$ is closed.
\end{Lemma}

{\it Proof.} Let $\TTT$ denote the restriction of $\TTT'$ to an automorphism of $\Cl$, which exists by assumption~(1).

Write $\Cl\hraa{i}\Cl'$ for the inclusion functor.

Since $\Cl'$ is closed, assumption (2) allows to conclude that $\Cl$ is a full additive subcategory of $\Cl'$, and moreover, that $\Cl$ is weakly abelian such that $i$ is subexact; cf.\ Lemma~\ref{LemC0_5}.

Given $n\ge 0$ and $X\in\Ob\,\ulk{\Cl^+(\b\De_n^\#)}$, we define, by restriction, $([X]^{+1}\lraa{X\tht_n} [X^{+1}]) := ([X]^{+1}\lrafl{28}{X\tht'_n} [X^{+1}])$. Since 
$\ulk{\Cl^+(\b\De_n^\#)}\;\lraa{\ul{i}}\;\ulk{\Cl'^+(\b\De_n^\#)}\ru{-2.3}$ is full and faithful by Remark~\ref{RemRes0}, this is a welldefined transformation satisfying $\tht_n\st\ul{i} = \ul{i}\st\tht'_n\,$.

Given $m,\, n\,\ge\, 0$ and a periodic monotone map $\b\De_n\llaa{p}\b\De_m\,$, we have $\ul{p}^\#\st\ul{i} = \ul{i}\st\ul{p}^\#$, whence 
\[
\ul{p}^\#\st\tht_m\st\ul{i} \= \ul{p}^\#\st\ul{i}\st\tht'_m \= \ul{i}\st\ul{p}^\#\st\tht'_m \= \ul{i}\st\tht'_n\st\ul{p}^\# \= \tht_n\st\ul{i}\st\ul{p}^\# \= \tht_n\st\ul{p}^\#\st\ul{i}\; ,
\]
so that we may conclude that $\ul{p}^\#\st\tht_m = \tht_n\st\ul{p}^\#$, for $\ul{i}$ is faithful.

Given $n\ge 0$, we have $\ul{\ffk}_n\!\st\ul{i} = \ul{i}\st\ul{\ffk}_n\,$, whence 
\[
\ul{\ffk}_n\st\tht_{n+1}\st\ul{i} \= \ul{\ffk}_n\st\ul{i}\st\tht'_{n+1} \= \ul{i}\st\ul{\ffk}_n\st\tht'_{n+1} \= \ul{i}\st\tht'_{2n+1}\st\ul{\ffk}_n 
\= \tht_{2n+1}\st\ul{i}\st\ul{\ffk}_n \= \tht_{2n+1}\st\ul{\ffk}_n\st\ul{i}\; ,
\]
so that we may conclude that $\ul{\ffk}_n\st\tht_{n+1} = \tht_{2n+1}\st\ul{\ffk}_n\,$, for $\ul{i}$ is faithful.

Hence $\tht$ is a Heller triangulation on $(\Cl,\TTT)$; cf.\ \bfcite{Ku05}{Def.\ 1.5.(i)}. By construction, $\Cl\hraa{i}\Cl'$ is strictly exact. 

By~(2) and Remark \ref{RemCl2}, the Heller triangulated category $(\Cl,\TTT,\tht)$ is closed.
\qed


\section{Functors are exact if and only if they are compatible with $n$-triangles}
\label{SecExCrit}

Suppose given Heller triangulated categories $(\Cl,\TTT,\tht)$ and $(\Cl',\TTT',\tht')$. 

Concerning the notion of \mb{$n$-triangles} in a Heller triangulated category, cf.\ \bfcite{Ku05}{Def.\ 1.5.(ii)}.

For $n\ge 0$, an object $Y$ in $\Cl(\bDes{n})$ is called {\it periodic} if $[Y]^{+1} = [Y^{+1}]$.

Suppose given an additive functor $\Cl\lraa{F}\Cl'$ and an isomorphism $\TTT F\lraisoa{a} F\TTT'$. 

For $z\in\Z$, we let $\TTT^z F\lraisofl{26}{a^{(z)}} F\TTT'^z$ be defined by 
\[
\ba{lcll}
a^{(0)}     & := & 1_F & \\
a^{(z + 1)} & := & (\TTT\st\, a^{(z)})\cdot (a\st\TTT'^z)                       & \text{for $z\ge 0$} \\
a^{(z - 1)} & := & (\TTT^-\st\, a^{(z)})\cdot (\TTT^-\st\, a^-\st\TTT'^{z-1})   & \text{for $z\le 0$} \\
\ea
\]
Then $(\TTT^z\st\,a^{(w)})\cdot (a^{(z)}\st\TTT'^w) = a^{(z+w)} : \TTT^{z+w} F\lraiso F\TTT'^{z+w}$ for $z,\, w\,\in\,\Z$.

Given a periodic $n$-pretriangle $X\in\Ob\,\Cl^{+,\,\per}(\bDes{n})$, for sake of brevity we denote in this section by 
\[
Y \;:=\; X(F(\bDes{n}))\;\in\;\Ob\,\Cl'(\bDes{n}) 
\]
the diagram obtained by pointwise application of $F$ to $X$. We have 
\[
Y|_{\dDe_n^{+1}} \= X|_{\dDe_n}((\TTT F)(\dDe_n)) \;\;\mraisofl{29}{X|_{\dDe_n}(a(\dDe_n))}\;\; X|_{\dDe_n} ((F\TTT')(\dDe_n)) \= (Y|_{\dDe_n})^{+1}\; .
\]
Isomorphic substitution along this isomorphism turns $Y|_{\b\De_n^{\tru\trd}}$ into a diagram $\br Y|_{{\b\De_n^{\tru\trd}}}$ for a periodic object $\br Y\in\Ob\,\Cl(\bDes{n})$ 
thus defined. We have an isomorphism $Y\lraisoa{\br a} \br Y$ in $\Cl'(\bDes{n})$ that at $(\be/\al)^{+z}$ for $0\le\al\le\be\le n$ and $z\in\Z$ is given by
\[
\Big(Y_{(\be/\al)^{+z}}\mraisofl{30}{\br a_{(\be/\al)^{+z}}} \br Y_{(\be/\al)^{+z}}\Big) 
\; :=\; \Big(X_{\be/\al} \TTT^z F\mraisofl{27}{X_{\be/\al}\, a^{(z)}} X_{\be/\al}F\TTT'^z\Big)\; .
\]
In fact, given $0\le\al\le n$ and $z\in\Z$, we obtain a commutative quadrangle
\[
\xymatrix{
X_{n/\al}\TTT^z F\ar[rrr]^{x\TTT^z F}\ar[d]_{X_{n/\al}\,a^{(z)}}^\wr & & & X_{\al/0} \TTT^{z+1} F\ar[d]^{X_{\al/0}\,a^{(z+1)}}_\wr \\
X_{n/\al}F\TTT'^z\ar[rrr]^{(xF\TTT'^z)(X_{\al/0}\,a \TTT'^z)}        & & & X_{\al/0} F\TTT'^{z+1}\zw{,}                            \\
}
\]
for
\[
(X_{n/\al}\,a^{(z)})(xF\TTT'^z)(X_{\al/0}\,a \TTT'^z)
\; =\; (x\TTT^z F)(X_{\al/0}\TTT a^{(z)})(X_{\al/0}\,a \TTT'^z) 
\; =\; (x\TTT^z F)(X_{\al/0}\,a^{(z+1)}) \; . 
\]

The remaining commutativities required for the naturality of $Y\lraisoa{\br a}\br Y$ follow by naturality of~$a^{(z)}$. 

We remark that $\br a|_{\dDe_n} = 1_{F(\dDe_n)}\,$.

If $F$ is subexact, then $Y$ is an $n$-pretriangle and $\br Y$ is a periodic $n$-pretriangle.

\begin{Lemma}
\label{LemC0_9}
Suppose given an exact functor $(F,a)$. 

Then for each $n$-triangle $X$ of $\Cl$, i.e.\ $X\in\Ob\Cl^{\tht=1,+}(\bDes{n})$, the object $\,\br Y$ of $\,\Cl'(\bDes{n})$ defined by {\rm (1)} and {\rm (2)} is an $n$-triangle of $\,\Cl'$, 
i.e.\ $\br Y\in\Ob\Cl'^{\,\tht=1,+}(\bDes{n})$. 

\begin{itemize}
\item[{\rm (1)}] We have $[\br Y]^{+1} = [\br Y^{+1}]$.
\item[{\rm (2)}] On $\b\De_n^{\tru\trd}$, the object $\br Y|_{{\b\De_n^{\tru\trd}}}$ arises from $Y := X(F(\bDes{n}))|_{{\b\De_n^{\tru\trd}}}$ by isomorphic substitution along
$Y|_{\dDe_n^{+1}} = X|_{\dDe_n}(\TTT(\dDe_n))(F(\dDe_n)) \;\;\lraisoa{a(\dDe_n)}\;\; X|_{\dDe_n} (F(\dDe_n))(\TTT'(\dDe_n)) = (Y|_{\dDe_n})^{+1}\ru{5.5}$.
\end{itemize}
\end{Lemma}

\bq
 Cf.\ \bfcite{Ku05}{Lem.\ 3.8} for the case of a strictly exact functor.
\eq

{\it Proof.}
Suppose given $n\ge 0$ and an $n$-triangle $X\in\Ob\,\Cl^{+,\,\tht = 1}(\bDes{n})$. By construction, $\br Y$ is periodic. We have to show
that $\br Y\tht'_n \sollgl 1_{[\br Y]^{+1}}\,$. We obtain
\[
\barcl
Y\tht'_n 
& = & X \left(\ulk{F^+(\bDes{n})}\st\tht'_n\right)                                       \vsp{2}\\
& = & X \left(\left(\tht_n\st\ulk{F^+(\bDes{n})}\right)\cdot \ulk{a^+(\bDes{n})}\right)  \vsp{2}\\
& = & X \left(\tht_n\st\ulk{F^+(\bDes{n})}\right)\cdot X\ulk{a^+(\bDes{n})}              \vsp{2}\\ 
& = & X\ulk{a^+(\bDes{n})} \; .                                                                 \\ 
\ea
\]
In particular, $Y\tht'|_{\dDe_n} = X|_{\dDe_n}\ulk{a(\dDe_n)}$. Hence, restricting the stably commutative quadrangle
\[
\xymatrix{
[Y]^{+1}\ar[r]^{[\br a]^{+1}}\ar[d]_{Y\tht'_n} & [\br Y]^{+1} \ar[d]^{\br Y\tht'_n} \\
[Y^{+1}]\ar[r]^{[\br a^{+1}]}                  & [\br Y^{+1}]                       \\
}
\]
to $\dDe_n\,$, we obtain the stably commutative quadrangle
\[
\xymatrix@C=15mm{
Y|_{\dDe_n^{+1}}\ar[r]^(0.43){X|_{\dDe_n}\ulk{a(\dDe_n)}}\ar[d]_{X|_{\dDe_n}\ulk{a(\dDe_n)}} & (Y|_{\dDe_n})^{+1} \ar[d]^{\br Y\tht'_n|_{\dDe_n}} \\
(Y|_{\dDe_n})^{+1}\ar[r]^{1}                                                                 & (Y|_{\dDe_n})^{+1}\zw{.}                           \\
}
\]
whence $\br Y\tht'_n|_{\dDe_n} = 1_{(Y|_{\dDe_n})^{+1}}\,$. Since the functor from $\ulk{\Cl^+(\bDes{n})}$ to $\ulk{\Cl(\dDe_n)}$ induced by restriction is an equivalence by 
\bfcite{Ku05}{Prop.\ 2.6}, this implies that $\br Y\tht'_n = 1_{[\br Y]^{+1}}\,$.
\qed

\begin{Proposition}
\label{PropC1}
Suppose $\Cl$ to be closed. 

The pair $(F,a)$ is an exact functor if and only if for each $n$-triangle $X$ of $\Cl$, the object $\,\br Y$ of $\,\Cl'(\bDes{n})$ defined by {\rm (1,\,2)} is an $n$-triangle of $\Cl'$. 
\begin{itemize}
\item[{\rm (1)}] We have $[\br Y]^{+1} = [\br Y^{+1}]$.
\item[{\rm (2)}] On $\b\De_n^{\tru\trd}$, the object $\br Y|_{{\b\De_n^{\tru\trd}}}$ arises from $Y := X(F(\bDes{n}))|_{{\b\De_n^{\tru\trd}}}$ by isomorphic substitution along
$Y|_{\dDe_n^{+1}} = X|_{\dDe_n}(\TTT(\dDe_n))(F(\dDe_n)) \;\;\lraisoa{a(\dDe_n)}\;\; X|_{\dDe_n} (F(\dDe_n))(\TTT'(\dDe_n)) = (Y|_{\dDe_n})^{+1}\ru{5.5}$.
\end{itemize}
\end{Proposition}

{\it Proof.} In view of Lemma~\ref{LemC0_9}, it suffices to show that if each $n$-triangle $X$ in $\Cl$ yields an $n$-triangle $\br Y$ in $\Cl'$ by (1,\,2), then $(F,a)$ is exact.

\hypertarget{claim.1}{}%
We {\it claim} that $F$ is subexact. By Lemma~\ref{LemC0_5}, it suffices to show that given a morphism $S\lraa{p} T$ in $\Cl$, there exists a weak cokernel of $p$ that is mapped by $F$
to a weak cokernel. Since $\Cl$ is a {\sf closed} Heller triangulated category, a weak cokernel of $p$ is contained in the the completion of $S\lraa{p} T$ to a $2$-triangle $X$ by Lemma~\ref{LemCl3}. 
We form the corresponding $2$-triangle $\br Y\,$ defined by~(1,\,2). Since it contains a weak cokernel of $SF\lraa{pF} TF$, and since $\br Y$ is isomorphic, in 
$\Cl'^+(\bDes{n})$, to $X(F^+(\b\De_2^\#))$, the image under $F$ of the weak cokernel of $p$ that is contained in the $2$-triangle $X$ is in fact a weak cokernel of $pF$. This proves the \hyperlink{claim.1}{\it claim.}

\hypertarget{claim.2}{}%
We {\it claim} that 
\[
(\tht_n\st \ulk{F^+(\bDes{n})})\cdot\ulk{a^+(\bDes{n})} \= \ulk{F^+(\bDes{n})}\st\tht'_n
\] 
for all $n\ge 0$. Suppose given $X\in\Ob\,\Cl^+(\bDes{n})$. Since $\Cl$ is a 
{\sf closed} Heller triangulated category, 
there exists an $n$-triangle $X'$ such that $X'|_{\dDe_n} = X|_{\dDe_n}\,$; cf.\ Lemma~\ref{LemCl3}. 
By \bfcite{Ku05}{Prop.\ 2.6}, we have an isomorphism $X\lraisoa{f} X'$ in $\ulk{\Cl^+(\bDes{n})}$ that restricts to the identity on $\dDe_n\,$. 
We dispose of a commutative diagram
\[
\xymatrix{
[X]^{+1}\ar[r]^{X\tht_n}\ar[d]_{[f]^{+1}} & [X^{+1}]\ar[d]^{[f^{+1}]} \\
[X']^{+1}\ar[r]^{X'\tht_n}                & [X'^{+1}]                 \\
}
\]
in $\ulk{\Cl^+(\bDes{n})}$. Since, by construction, $X'\tht_n = 1$, we have $X\tht_n = [f]^{+1}\cdot[f^{+1}]^-$ in $\ulk{\Cl^+(\bDes{n})}$. 

Likewise, we have a commutative quadrangle
\[
\xymatrix@C=20mm{
[X\ulk{F^+(\b\De_n)}]^{+1}\ar[r]^(0.47){X\ulk{F^+(\b\De_n)}\tht'_n}\ar[d]_{[f\ulk{F^+(\b\De_n)}]^{+1}} & [(X\ulk{F^+(\b\De_n)})^{+1}]\ar[d]^{[(f\ulk{F^+(\b\De_n)})^{+1}]} \\
[X'\ulk{F^+(\b\De_n)}]^{+1}\ar[r]^(0.47){X'\ulk{F^+(\b\De_n)}\tht'_n}                                  & [(X'\ulk{F^+(\b\De_n)})^{+1}] \; ,                                \\
}
\]
in $\ulk{\Cl'^+(\bDes{n})}$. We want to calculate its lower arrow. Since $X'$ is an $n$-triangle, we have an isomorphism $Y'\lraisoa{\br a'} \br Y'$ formed as above, where $\br Y'\tht'_n = 1$. 
The stably commutative quadrangle
\[
\xymatrix{
[Y']^{+1}\ar[r]^{[\br a']^{+1}}\ar[d]_{Y'\tht'_n} & [\br Y']^{+1} \ar[d]^{\br Y'\tht'_n \;=\; 1} \\
[Y'^{+1}]\ar[r]^{[\br a'^{+1}]}                   & [\br Y'^{+1}]                               \\
}
\]
yields by restriction to $\dDe_n$ the commutative diagram
\[
\xymatrix@C+15mm{
X'|_{\dDe_n}\,\ulk{(\TTT F)(\dDe_n)}\ar[r]^-{X'|_{\dDe_n}\,\ulk{a(\dDe_n)}}\ar@{=}[d] & X'|_{\dDe_n}\,\ulk{(F\,\TTT')(\dDe_n)}\ar@{=}[d] \\
[Y']^{+1}|_{\dDe_n}\ar[r]^{[\br a]^{+1}|_{\dDe_n}}\ar[d]_{Y'\tht'_n|_{\dDe_n}}        & [\br Y']^{+1}|_{\dDe_n}\ar[d]^1                  \\
[Y'^{+1}]|_{\dDe_n}\ar[r]^1                                                           & [\br Y'^{+1}]|_{\dDe_n} \zw{,}                   \\ 
}
\]
whence $Y'\tht'_n|_{\dDe_n} = X'|_{\dDe_n}\,\ulk{a(\dDe_n)} = X'\ulk{a^+(\bDes{n})}|_{\dDe_n}\,$. Since the functor from $\ulk{\Cl'^+(\bDes{n})}$ to $\ulk{\Cl'(\dDe_n)}$ induced by 
restriction is an equivalence by \bfcite{Ku05}{Prop.\ 2.6}, this implies that 
\[
X'\ulk{F^+(\b\De_n)}\tht'_n \= Y'\tht'_n \= X'\ulk{a^+(\bDes{n})} \; .
\]
So we can conclude that
\[
\barcl
X\ulk{F^+(\b\De_n)}\tht'_n 
& = & [f\ulk{F^+(\b\De_n)}]^{+1}\cdot X'\ulk{F^+(\b\De_n)}\tht'_n\cdot[(f\ulk{F^+(\b\De_n)})^{+1}]^-  \vsp{2} \\
& = & [f\ulk{F^+(\b\De_n)}]^{+1}\cdot X'\ulk{a^+(\bDes{n})}\cdot[(f\ulk{F^+(\b\De_n)})^{+1}]^-        \vsp{2} \\
& = & [f\ulk{F^+(\bDes{n})}]^{+1}\cdot X'\ulk{a^+(\bDes{n})} \cdot (f\ulk{(F \TTT')^+(\bDes{n})})^-   \vsp{2} \\ 
& = & [f\ulk{F^+(\bDes{n})}]^{+1}\cdot (f\ulk{(\TTT F)^+(\bDes{n})})^-\cdot X\ulk{a^+(\bDes{n})}      \vsp{2} \\ 
& = & ([f]^{+1}\cdot [f^{+1}]^-)\ulk{F^+(\bDes{n})}\cdot X\ulk{a^+(\bDes{n})}                         \vsp{2} \\ 
& = & X(\tht_n\st \ulk{F^+(\bDes{n})})\cdot X\ulk{a^+(\bDes{n})} \; .                                         \\ 
\ea
\]
This proves the \hyperlink{claim.2}{\it claim.}
\qed

\begin{Corollary}
\label{CorC2}
Suppose $(\Cl,\TTT,\tht)$ to be closed.

Suppose given an additive functor $\Cl\lraa{F}\Cl'$ such that $\TTT F = F\,\TTT'$.

Then $F$ is strictly exact if and only if for each $n\ge 0$ and each $n$-triangle $X\in\Ob\,\Cl^{+,\,\tht = 1}(\bDes{n})$, the diagram $X(F(\bDes{n}))\in\Ob\,\Cl'(\bDes{n})$, obtained
by pointwise application of $F$, is an $n$-triangle.
\end{Corollary}

{\it Proof.} In this case, we have $a = 1_{\TTT F} = 1_{F\,\TTT'}$ and $\br Y = Y = X(F(\bDes{n}))$. Since $(F,1)$ is exact if and only if $F$ is strictly exact, the assertion follows by 
Proposition \ref{PropC1}.\qed


\section{Adjoints}
\label{SecAdj}

\subsection{Adjoints and shifts}
\label{SecAdjSh}

Suppose given categories $\Al$ and $\Al'$. Suppose given an endofunctor $T$ of $\Al$. Suppose given an endofunctor $T'$ of $\Al'$.

Suppose given functors $\Al\lradoubleinva{F}{G}\Al'\ru{-5}$ such that $F\adj G$ via unit $1\lraa{\eps} FG$ and counit $GF\lraa{\et} 1$, i.e.\ $(G\eps)(\et G) = 1_G$ and $(\eps F)(F\et) = 1_F\,$.

Suppose given $TF\lraa{\al} FT'$. 

Let
\[
(GT\lraa{\be} T'G) \; :=\; \big(GT \lrafl{28}{GT\eps} GTFG \lrafl{28}{G\al G} GFT'G\lrafl{30}{\et T'G} T'G\big) \; .
\]
So we have the commutative diagram
\[
\xymatrix@C+5mm{
GT\ar[r]^\be\ar[d]_{GT\eps} & T'G                         \\
GTFG\ar[r]^{G\al G}         & GFT'G\ar[u]_{\et T'G}\zw{.} \\ 
}
\]

\pagebreak

\begin{Lemma}
\label{LemAdPrep}\Absit
\begin{itemize}
\item[\rm (1)] We have the commutative diagram
\[
\xymatrix@C+5mm{
TF\ar[r]^\al\ar[d]_{\eps TF} & FT'                        \\
FGTF\ar[r]^{F\be F}          & FT'GF\ar[u]_{FT'\et}\zw{.} \\ 
}
\]
\item[\rm (2)] We have the commutative quadrangle
\[
\xymatrix{
T\ar[r]^{\eps T}\ar[d]_{T\eps} & FGT\ar[d]^{F\be} \\
TFG\ar[r]_{\al G}              & FT'G \zw{.}      \\
}
\]
\item[\rm (2$^\0$)] We have the commutative quadrangle
\[
\xymatrix{
GTF \ar[r]^{G\al}\ar[d]_{\be F} & GFT'\ar[d]^{\et T'} \\
T'GF\ar[r]_{T'\et}              & T' \zw{.}           \\
}
\]
\item[\rm (3)] Suppose that $T$ and $T'$ are autofunctors. Write $G' = T'GT^-$. If $\al$ is an isotransformation, then so is $\be$, where 
\[
\ba{l}
(T'G \lrafl{28}{\be^-} GT) \= \\
\big( T'G = T'GT^- T \mrafl{28}{T'GT^-\eps T} T'GT^-FGT = \\ 
T'GT^-FT'T'^-GT \;\;\mrafl{28}{T'GT^-\al^-T'^-GT}\;\; T'GT^-TFT'^-GT = \\
T'GFT'^-GT \mrafl{28}{T'\et T'^-GT} T'T'^-GT = GT\big)\; .
\ea
\]

\item[\rm (3$^\0$)] Suppose that $T$ and $T'$ are autofunctors. If $\be$ is an isotransformation, then so is $\al$, where 
\[
\ba{l}
(FT' \lrafl{28}{\al^-} TF) \= \\
\big( FT' = TT^-FT' \mrafl{28}{T\eps T^- FT'} TFGT^-FT' = \\ 
TFT'^-T'GT^-FT' \;\;\mrafl{28}{TFT'^-\be^-T^-FT'}\;\; TFT'^-GTT^-FT' = \\
TFT'^-GFT' \mrafl{28}{TFT'^-\et T'} TFT'^-T' = TF \big)\; .
\ea
\]
\end{itemize}
\end{Lemma}

{\it Proof.} Ad (2). We have
\[
\barcl
(\eps T)(F\be)
& = & (\eps T)(FGT\eps)(FG\al G)(F\et T'G) \\
& = & (T\eps)(\eps TFG)(FG\al G)(F\et T'G) \\
& = & (T\eps)(\al G)(\eps FT'G)(F\et T'G)  \\
& = & (T\eps)(\al G) \; .                  \\
\ea
\]

Ad (1). We have
\[
\barcl
(\eps TF)(F\be F)(FT'\et)
& \aufgl{(2)} & (T\eps F)(\al GF)(FT'\et) \\
& =           & (T\eps F)(TF\et) \al      \\
& =           & \al\; .                   \\
\ea
\]
\qed

Ad (3). We have
\[
\ba{cl}
  & \be\cdot (T'GT^-\eps T)(T'GT^-\al^-T'^-GT)(T'\et T'^-GT)          \\
= & (\be T^-T)(T'GT^-\eps T)(T'GT^-\al^-T'^-GT)(T'\et T'^-GT)         \\
= & (GTT^-\eps T)(\be T^-FGT)(T'GT^-\al^-T'^-GT)(T'\et T'^-GT)        \\
= & (G\eps T)(\be T^-FT'T'^-GT)(T'GT^-\al^-T'^-GT)(T'\et T'^-GT)      \\
= & (G\eps T)(GTT^-\al^-T'^-GT)(\be T^-TFT'^-GT)(T'\et T'^-GT)        \\
= & (G\eps T)(G\al^-T'^-GT)(\be FT'^-GT)(T'\et T'^-GT)                \\
\aufgl{(2$^\0$)} & (G\eps T)(G\al^-T'^-GT)(G\al T'^-GT)(\et T'T'^-GT) \\
= & (G\eps T)(\et GT)                                                 \\
= & 1                                                                 \\
\ea
\]
and
\[
\ba{cl}
  & (T'GT^-\eps T)(T'GT^-\al^-T'^-GT)(T'\et T'^-GT)\cdot\be       \\
= & (T'GT^-\eps T)(T'GT^-\al^-T'^-GT)(T'\et T'^-GT)(T'T'^-\be)    \\
= & (T'GT^-\eps T)(T'GT^-\al^-T'^-GT)(T'GFT'^-\be)(T'\et T'^-T'G) \\
= & (T'GT^-\eps T)(T'GT^-\al^-T'^-GT)(T'GT^-TFT'^-\be)(T'\et G)   \\
= & (T'GT^-\eps T)(T'GT^-FT'T'^-\be)(T'GT^-\al^-T'^-T'G)(T'\et G) \\
= & (T'GT^-\eps T)(T'GT^-F\be)(T'GT^-\al^-G)(T'\et G)             \\
\aufgl{(2)} & (T'GT^-T\eps)(T'GT^-\al G)(T'GT^-\al^-G)(T'\et G)   \\
= & (T'G\eps)(T'\et G)                                            \\
= & 1 \; .                                                        \\
\ea
\]

\kommentar{
Ad (3^\0). We have
\[
\ba{cl}
  & \al\cdot (T\eps T^- FT')(TFT'^-\be^-T^-FT')(TFT'^-\et T')        \\
= & (TT^-\al)(T\eps T^- FT')(TFT'^-\be^-T^-FT')(TFT'^-\et T')        \\
= & (T\eps T^-TF)(TFGT^-\al)(TFT'^-\be^-T^-FT')(TFT'^-\et T')        \\
= & (T\eps F)(TFT'^-T'GT^-\al)(TFT'^-\be^-T^-FT')(TFT'^-\et T')      \\
= & (T\eps F)(TFT'^-\be^-T^-TF)(TFT'^-GTT^-\al)(TFT'^-\et T')        \\
= & (T\eps F)(TFT'^-\be^-F)(TFT'^-G\al)(TFT'^-\et T')                \\
\aufgl{(2$^\0$)} & (T\eps F)(TFT'^-\be^-F)(TFT'^-\be F)(TFT'^-T'\et) \\
= & (T\eps F)(TF\et)                                                 \\
= & 1                                                                \\
\ea
\]
and
\[
\ba{cl}
  & (T\eps T^- FT')(TFT'^-\be^-T^-FT')(TFT'^-\et T')\cdot\al        \\
= & (T\eps T^- FT')(TFT'^-\be^-T^-FT')(TFT'^-\et T')(\al T'^-T')    \\
= & (T\eps T^- FT')(TFT'^-\be^-T^-FT')(\al T'^-GFT')(FT'T'^-\et T') \\
= & (T\eps T^- FT')(TFT'^-\be^-T^-FT')(\al T'^-GTT^-FT')(F\et T')   \\
= & (T\eps T^- FT')(\al T'^-T'GT^-FT')(FT'T'^-\be^-T^-FT')(F\et T') \\
= & (T\eps T^- FT')(\al GT^-FT')(F\be^-T^-FT')(F\et T')             \\
\aufgl{(2)} & (\eps TT^- FT')(F\be T^-FT')(F\be^-T^-FT')(F\et T')   \\
= & (\eps FT')(F\et T')                                             \\
= & 1\; .                                                           \\
\ea
\]
}

\subsection{An adjoint of an exact functor is exact}
\label{SecAdjEx}

The Verdier triangulated version of the following proposition is due to {\sc Margolis} \mb{\bfcite{Ma83}{App.\ 2, Prop.\ 11}}, and, in a more general form, to {\sc Keller} and {\sc Vossieck} 
\bfcite{KV87}{1.6}.

\begin{Proposition}
\label{PropAd1}
Suppose given Heller triangulated categories $(\Cl,\TTT,\tht)$ and $(\Cl',\TTT',\tht')$. 

Suppose given an exact functor $(F,a)$ from $\Cl$ to $\Cl'\,$; cf.\ Definition~{\rm\ref{DefExact}.} 

Suppose given a functor $\,\Cl\llaa{G}\Cl'\,$. 

So $\;\Cl\lradoubleinva{F}{G}\Cl'\;$ and $\;\TTT F \lraisoa{a} F\TTT'\;$.

\begin{itemize}
\item[\rm (1)] If $F\adj G$, then there exists an isomorphism $\TTT' G \lraisoa{b} G\TTT$ such that $(G,b)$ is an exact functor from $\Cl'$ to $\Cl$.

Choose a unit $1_{\Cl'}\lraa{\eps} FG$ and a counit $GF\lraa{\eta} 1_{\Cl}\,$. Then, more precisely, we may choose 
\[
(\TTT' G\;\lraa{b}\; G\TTT) \;\; :=\;\; (G\TTT\;\lrafl{27}{G\TTT\eps}\; G\TTT FG \;\lrafl{27}{G a G}\; GF\TTT' G \;\lrafl{27}{\et\TTT' G}\; \TTT'G)^-\; .
\]

\item[\rm (1$^\0$)] If $G\adj F$, then there exists an isomorphism $\TTT' G \lraisoa{b} G\TTT$ such that $(G,b)$ is an exact functor from $\Cl'$ to $\Cl$.

Choose a unit $1_{\Cl'}\lraa{\eps} GF$ and a counit $FG\lraa{\eta} 1_{\Cl}\,$. Then, more precisely, we may choose 
\[
(\TTT' G\;\lraa{b}\; G\TTT) \;\; :=\;\; (\TTT' G\;\lrafl{27}{\eps\TTT'G}\; GF\TTT' G \;\lrafl{27}{G a^- G}\; G\TTT FG \;\lrafl{27}{G\TTT\eta}\; G\TTT)\; .
\]
\end{itemize}
\end{Proposition}

{\it Proof.} Ad (1). By Lemma~\ref{LemAdjSubex}.(1$^\0$), $G$ is subexact.

Lemma~\ref{LemAdPrep}.(3) yields the isotransformation $b^- := (G\TTT\eps)(GaG)(\et\TTT' G)$.

Suppose given $n\ge 0$. We shall make use of the abbreviation $\ulk{G} = \ulk{G^+(\bDes{n})}\,$, etc. We have to show
that 
\[
(\tht'_n\st\ulk{G})\cdot\ulk{b} \;\sollgl\; \ulk{G}\st\tht_n\; ,
\]
i.e.\ that
\[
(\ulk{G}\st\tht_n)\cdot\ulk{b}^- \;\sollgl\; \tht'_n\st\ulk{G} \; ,
\]
i.e.\ that
\[
(\ulk{G}\st\tht_n)\cdot (\ulk{G}\st\ulk{\TTT}\st\ulk{\eps})\cdot (\ulk{G}\st \ulk{a} \st \ulk{G})\cdot (\ulk{\et}\st\ulk{\TTT'}\st\ulk{G}) \;\sollgl\; \tht'_n\st\ulk{G}
\]

Recall that $[-]^{+1}$ denotes the outer shift, that $[-^{+1}]$ denotes the inner shift and that \linebreak 
$\tht_n : [-]^{+1} \lraiso [-^{+1}]$ on $\ulk{\Cl^+(\bDes{n})}\,$; similarly on $\ulk{\Cl'^+(\bDes{n})}\,$.

We obtain
\[
\ba{cl}
  & (\ulk{G}\st\tht_n)\cdot (\ulk{G}\st\ulk{\TTT}\st\ulk{\eps})\cdot (\ulk{G}\st \ulk{a} \st \ulk{G})\cdot (\ulk{\et}\st\ulk{\TTT'}\st\ulk{G})                      \\
= & (\ulk{G}\st\tht_n)\cdot (\ulk{G}\st [-^{+1}]\st\ulk{\eps})\cdot (\ulk{G}\st \ulk{a} \st \ulk{G})\cdot (\ulk{\et}\st\ulk{\TTT'}\st\ulk{G})                       \\
= & (\ulk{G}\st [-]^{+1}\st\ulk{\eps})\cdot (\ulk{G}\st\tht_n\st \ulk{F}\st \ulk{G})\cdot (\ulk{G}\st \ulk{a} \st \ulk{G})\cdot (\ulk{\et}\st\ulk{\TTT'}\st\ulk{G}) \\
= & (\ulk{G}\st [-]^{+1}\st\ulk{\eps})\cdot (\ulk{G}\st ((\tht_n\st \ulk{F})\cdot\ulk{a})\st \ulk{G})\cdot (\ulk{\et}\st\ulk{\TTT'}\st\ulk{G})                      \\
\aufgl{$(F,a)$ ex.} & (\ulk{G}\st [-]^{+1}\st\ulk{\eps})\cdot (\ulk{G}\st \ulk{F}\st\tht'_n\st \ulk{G})\cdot (\ulk{\et}\st\ulk{\TTT'}\st\ulk{G})                    \\
= & (\ulk{G}\st [-]^{+1}\st\ulk{\eps})\cdot (\ulk{G}\st \ulk{F}\st\tht'_n\st \ulk{G})\cdot (\ulk{\et}\st [-^{+1}]\st\ulk{G})                                        \\
= & (\ulk{G}\st [-]^{+1}\st\ulk{\eps})\cdot (\ulk{\et}\st [-]^{+1}\st\ulk{G})\cdot (\tht'_n\st \ulk{G})                                                             \\
= & (\ulk{G}\st\ulk{\eps}\st [-]^{+1})\cdot (\ulk{\et}\st\ulk{G}\st [-]^{+1})\cdot (\tht'_n\st \ulk{G})                                                             \\
= & \tht'_n\st \ulk{G}\; .                                                                                                                                          \\
\ea
\]

Ad (1$^\0$). Cf.\ Lemma~\ref{LemAdPrep}.(1). 
\qed

\begin{Example}
\label{ExAd1_5}\rm
Suppose we are in the situation of Proposition~\ref{PropAd1}.(1). Then $\eps$ and $\et$ are periodic; cf.\ Definition~\ref{DefExactPeriodic}.

Ad $\eps : 1\lra FG$. The functor $(F,a)(G,b) = (FG\,,\,aG\cdot Fb)$ is exact; cf.\ Remark~\ref{RemD1}. The functor $(1_\Cl\,,\,1)$ is exact. The quadrangle
\[
\xymatrix{
\TTT\ar[r]^-{T\eps}\ar[d]_{1} & \TTT FG\ar[d]^{aG\cdot Fb} \\
\TTT\ar[r]^-{\eps T}          & FG\TTT                     \\ 
}
\]
commutes by Lemma~\ref{LemAdPrep}.(2).

Ad $\et : GF\lra 1$. The functor $(G,b)(F,a) = (GF\,,\,bF\cdot Ga)$ is exact; cf.\ Remark~\ref{RemD1}. The functor $(1_{\Cl'}\,,\,1)$ is exact. The quadrangle
\[
\xymatrix{
\TTT' GF\ar[r]^-{\TTT'\et}\ar[d]_{bF\cdot Ga} & \TTT'\ar[d]^{1} \\
GF\TTT' \ar[r]^-{\et\TTT'}                    & \TTT'           \\
}
\]
commutes by Lemma~\ref{LemAdPrep}.(2$^\0$).
\end{Example}

\subsection{A functor shiftcompatibly adjoint to a strictly exact functor is strictly exact}
\label{SecAdjExStr}

Suppose given closed Heller triangulated categories $(\Cl,\TTT,\tht)$ and $(\Cl',\TTT',\tht')$ 

Recall that an additive functor $F : \Cl\lra\Cl'$ is strictly exact if and only if $(F,1)$ is exact; cf.\ \mb{\bfcite{Ku05}{Def.\ 1.5.(iii)}}, Definition \ref{DefExact}.

\begin{Corollary}
\label{CorAd2}\Absit

Suppose given a strictly exact functor $\,\Cl\lraa{F}\Cl'\,$. 

Suppose given a functor $\,\Cl\llaa{G}\Cl'\,$.

\begin{itemize}
\item[\rm (1)] If $F\adj G$, with unit $\eps : 1 \lra FG$ and counit $\et : GF \lra 1$ such that $(G\TTT\eps)(\et\TTT' G) = 1$, then $G$ is strictly exact.

\item[\rm (1$^\0$)] If $G\adj F$, with unit $\eps : 1 \lra FG$ and counit $\et : GF \lra 1$ such that $(\eps\TTT'G)(G\TTT\et) = 1$, then $G$ is strictly exact.
\end{itemize}
\end{Corollary}

{\it Proof.} Ad (1). In the notation of Proposition \ref{PropAd1}.(1), we have $a = 1$, and, consequently, 
$b = (G\TTT\eps)(\et\TTT' G) = 1$. Hence by loc.\ cit., $(G,1)$ is exact, i.e.\ $G$ is strictly exact.\qed


\section{Localisation}
\label{SecLoc}

We prove that the localisation $\Cl\fbby\Nl$ of a Heller triangulated category $\Cl$ at a thick subcategory $\Nl$ is Heller triangulated in such a way that the localisation functor 
$\Cl\lraa{\LLL}\Cl\fbby\Nl$ is strictly exact; cf.\ \bfcite{Ku05}{Def.\ 1.5}. There is considerable overlap 
with the classical localisation theory of Verdier triangulated categories, due to {\sc Verdier} \bfcit{Ve63}, which we include for sake of self-containedness.

Let $(\Cl,\TTT,\tht)$ be a closed Heller triangulated category; cf.\ Definition \ref{DefClosed}.

\begin{Definition}
\label{DefL0}\rm
A full additive subcategory $\Nl\tm\Cl$ is called {\it thick} if the conditions (1,\,2,\,3) are satisfied; cf.\ \bfcite{Ri89}{Prop.\ 1.3}
\begin{itemize}
\item[(1)] We have $\Nl^{+1} = \Nl$ (closed under shift).
\item[(2)] Given a $2$-triangle $(X,Y,Z)$ in $\Cl$ with $X$ and $Y$ in $\Ob\Nl$, then $Z\in\Ob\Nl$ \\ 
 (closed under taking cones). 
\item[(3)] Given $X,\, Y\in\Ob\Cl$ with $X\ds Y$ in $\Ob\Nl$, then $X\in\Ob\Nl$ \\
 (closed under taking summands).
\end{itemize}

Let $\Nl$ be a thick subcategory of $\Cl$. By Lemma \ref{LemS2}, conditions (1) and (2) of Definition \ref{DefL0} yield that $\Nl$ is a Heller triangulated subcategory of $\Cl$.

Let $M(\Nl) := \{ (X\lraa{f} Y)\in\Cl\; :\; \text{the cone of $f$ is in $\Ob\Nl$}\}$. An element of $M(\Nl)$ is called an {\it $M(\Nl)$-isomorphism} or often just an {\it \Niso}\ (not to 
be confused with ``an isomorphism in $\Nl\,$''). If $\Nl$ is unambiguous, then an \Niso\ is denoted by $X\imp Y$. For instance, $X\imp 0$ if and only if $0\imp X$ if and only if $X\in\Ob\Nl$.
\end{Definition}

\begin{Lemma}
\label{LemL1}
The subset $M(\Nl)$ of \Niso s in $\Cl$ is a multiplicative system in $\Cl$ in the sense of Definition~{\rm\ref{DefAppMult}.}
\end{Lemma}

{\it Proof.} Ad (Fr 2). Suppose given $X_{1/0}\lraa{x} X_{2/0}\lraa{x} X_{3/0}\lraa{x} X_{4/0}$ such that $X_{1/0}\lraa{x} X_{3/0}$ and $X_{2/0}\lraa{x} X_{4/0}$ are $M$-isomorphisms. We complete
to a $4$-triangle $X\in\Ob\Cl^{+,\,\tht=1}(\b\De_n^\#)$ using closedness of $\Cl\,$; cf.\ Lemma~\ref{LemCl3}. By \mb{\bfcite{Ku05}{Lem.\ 3.4.(1,\,6)}}, we have 
\mb{$X_{3/1}\,,\, X_{4/2}\,\in\,\Ob\Nl$}. We have to show that \mb{$X_{2/1}\,,\, X_{3/2}\,,\, X_{4/3}\,,\, X_{4/1}\,\in\,\Ob\Nl$.} Let the periodic monotone map 
$\b\De_5\lraa{p}\b\De_4$ be defined by $0p := 1$, $1p := 1$, $2p := 2$, $3p := 3$, $4p := 4$ and $5p := 4$. The $2$-triangle $Xp^\#\ffk_2\in\Ob\Cl^{+,\,\tht=1}(\b\De^\#_2)$ is given by
\[
\hspace*{-10mm}
\xymatrix@R=11mm{
& & & & 0 \\
& & & 0\ar[r] &  X_{2^{+1}/4}\ar[u] \\
& & 0\ar[r] &  X_{1^{+1}/2}\ar[r]^{-x}\ar[u]\ar@{}[ur]|{+} & X_{1^{+1}/4}\ar[u]^{x} \\
& 0\ar[r] & X_{4/3}\ar[r]^(0.33){\smatez{1}{0}}\ar[u]\ar@{}[ur]|{+}
        & X_{4/3}\dk X_{1^{+1}/2}\ar[r]^(0.58){\rsmatze{x}{-x}}\ar[u]^{\smatze{0}{1}}\ar@{}[ur]|{+} 
        & X_{1^{+1}/3}\ar[u]^{x}  \\ 
0\ar[r] & X_{3/1}\ar[r]^{x}\ar[u]\ar@{}[ur]|{+} & X_{4/1}\ar[r]^{x}\ar[u]^{x}\ar@{}[ur]|{+}
        & X_{4/2}\ar[r]\ar[u]^{\smatez{x}{x}}\ar@{}[ur]|{+} & 0\ar[u]\zw{,} \\
}
\]
cf.\ \bfcite{Ku05}{Lem.\ 3.4.(1,\,2), \S 1.2.1.2, \S 1.2.2.2.}. Since $X_{3/1}\,,\, X_{4/2}\,\in\,\Ob\Nl$, and since $\Nl$ is closed under cones, we have $X_{4/3}\ds X_{1^{+1}/2}\in\Ob\Nl$.
Since $\Nl$ is {\sf closed under summands} and under shift, we obtain $X_{4/3}\,,\, X_{2/1}\,\in\,\Ob\Nl$. Since $\Nl$ is closed under cones and under shift, $X_{4/1}\in\Ob\Nl$ ensues. Considering
$X$ again, since $\Nl$ is closed under cones, we finally obtain $X_{3/2}\in\Ob\Nl$.

Ad (Fr 3). Let $X\lraa{f} Y$ be a morphism in $\Cl$ such that there exists \mb{$Y\auf{s}{\imp} Z$} with \mbox{$fs = 0$.} We obtain a factorisation $(X\lraa{f} Y) = (X\lraa{u} N\lraa{v} Y)$
with $N\in\Ob\Nl$. Com\-pleting $(T_{1/0}\lraa{t} T_{2/0}\lraa{t} T_{3/0}) := (X\lraa{u} N\lraa{v} Y)$ to a $3$-triangle by Lemma~\ref{LemCl3}, we obtain
\mb{$T_{1/2^{-1}}\auf{t}{\imp} T_{1/0}\,$}, which composes to zero with $(T_{1/0}\lraa{t} T_{3/0}) \= (X\lraa{f} Y)$.

Ad (Fr 4).  Suppose given 
\[
\xymatrix{
X'                       &   \\
X\ar@{=>}[r]\ar[u]       & Y \\
}
\]
in $\Cl$. Prolonging $X\lra X'$ to a $2$-triangle $(X'',X,X')$, then completing $X''\lra X\imp Y$ to a \mb{$3$-triangle} using Lemma~\ref{LemCl3}, we obtain, by \bfcite{Ku05}{Lem.\ 3.4.(6)}, a
\mb{$3$-triangle} $T$ with $(T_{2/0}\lraa{t} T_{3/0}) = (X\imp Y)$ and $(T_{2/0}\lraa{t} T_{2/1}) = (X\lra X')$. Then $T_{3/2}\in\Ob\Nl$, whence $T_{2/1}\imp T_{3/1}\,$. The weak square
$(T_{2/0}\,,T_{3/0}\,,T_{2/1}\,,T_{3/1})$ is a completion as sought.\qed

\bq
 Note that if $X\in(\Ob\Cl)\ohne (\Ob\Nl)$, then $(0,0,0,X)$ is a weak square in which $0\lra 0$ is an \mb{\Niso}, but $0\lra X$ is not.
\eq

The localisation of $\Cl$ at $M(\Nl)$, defined as in \S\ref{SubsecMultSyst}, is also called the {\it localisation of $\Cl$ at $\Nl$,} and also written $\CNl := \Cl_{M(\Nl)}$.
Concerning the {\it localisation functor} $\Cl\lraa{\LLL}\CNl$, we refer to \S\ref{SubsecMultSyst}.

Recall that an additive functor between weakly abelian categories is called subexact if it induces an exact functor on the Freyd categories; cf.\ \bfcite{Ku05}{\S1.2.1.3}; cf.\ also Lemma~\ref{LemC0_5}.

\begin{Lemma}
\label{LemL2}
The category $\CNl$ is weakly abelian. The functor $\Cl\lraa{\LLL}\CNl$ is subexact.
\end{Lemma}

{\it Proof.} By Remark \ref{RemA2}, the category $\CNl$ is additive, and the localisation functor $\LLL : \Cl\lra\CNl$ is additive. 
We claim that $\LLL$ maps weak kernels to weak kernels. Let $X\lraa{f} Y$ be a weak kernel of $Y\lraa{g} Z$ in 
$\Cl$. We claim that it remains a weak kernel in $\CNl$. Suppose given a morphism $T\lraa{t} Y$ in $\Cl$ such that $tg = 0$ in $\CNl$, which we, by isomorphic replacement, may assume given. 
Let $T'\auf{s}{\imp} T$ be such that $stg = 0$ in $\Cl\,$; cf.\ Remark~\ref{RemA0_1}. Since $f$ is 
a weak kernel of $g$ in $\Cl$, we have a factorisation $st = uf$. Hence $t = (s^{-1} u) f$ is a factorisation of $t$ over $f$ in $\CNl$. 

Substituting isomorphically in $\CNl$ and using duality, for $\CNl$ to be weakly abelian, it suffices to show that each morphism $X\lraa{f} Y$ has a weak kernel resp.\ is a weak kernel in $\CNl$. 
But by the property of $\LLL$ just shown, we may use a weak kernel of $f$ in $\Cl$ resp.\ a morphism $f$ is a weak kernel of in $\Cl$.\qed

\begin{Remark}
\label{RemL3}
The category $\CNl$ carries a shift automorphism $\CNl\lraa{\TTT}\CNl$, $f/t\lramaps f^{+1}/t^{+1}$. We have $\LLL\TTT = \TTT\LLL$.
\end{Remark}

{\it Proof.} This functor is welldefined since $\Nl$, and hence $M(\Nl)$, is closed under shift in $\Cl$. Likewise, its inverse $f/t\lramaps f^{-1}/t^{-1}$ is welldefined.\qed

\begin{Lemma}
\label{LemL4}
Suppose given a Heller triangulated category $(\Dl,\TTT,\theta)$.

Suppose given a weakly abelian category $\Dl'$ and an automorphism $\Dl'\lraa{\TTT'}\Dl'$. Suppose given a subexact additive functor $\Dl\lraa{G}\Dl'\ru{5.0}$ strictly compatible with shift, i.e.\ $G\TTT' = \TTT G$.
Suppose that $\,\Dl(\dDe_n)\lrafl{25}{G(\dDe_n)}\Dl'(\dDe_n)$ is $1$-epimorphic for $n\ge 0$. 

Then the functor $\,\ulk{\Dl^+(\b\De_n^\#)}\mrafl{32}{\ulk{G^+(\b\De_n^\#)}} \ulk{\Dl'^+(\b\De_n^\#)}$ is $1$-epimorphic.

Moreover, there exists a unique Heller triangulation $\theta'$ on $(\Dl',\TTT')$ such that $\Dl\lraa{G}\Dl'$ is strictly exact; cf.\ \mb{\bfcite{Ku05}{Def.\ 1.5}}.
\end{Lemma}

{\it Proof.} Given $n\ge 0$. Since the residue class functors $\Dl(\dDe_n)\lra\ulk{\Dl(\dDe_n)}$ and 
$\Dl'(\dDe_n)\lra\ulk{\Dl'(\dDe_n)}$ are full and dense, they are $1$-epimorphic by \mb{\bfcite{Ku05}{Cor.\ A.37}}\,; concerning notation, cf.\ \mb{\bfcite{Ku05}{\S2.4}}. The commutative quadrangle
\[
\xymatrix@C=13mm{
\Dl(\dDe_n)\ar[r]^{G(\dDe_n)}\ar[d]       & \Dl'(\dDe_n)\ar[d] \\
\ulk{\Dl(\dDe_n)}\ar[r]^{\ulk{G(\dDe_n)}} & \ulk{\Dl'(\dDe_n)} \\
}
\]
shows that $\ulk{G(\dDe_n)}$ is $1$-epimorphic. Restriction induces equivalences $\ulk{\Dl^+(\b\De_n^\#)}\;\lrafl{32}{(-)|_{\dDe_n}}\;\ulk{\Dl(\dDe_n)}$ and 
$\ulk{\Dl'^+(\b\De_n^\#)}\;\lrafl{32}{(-)|_{\dDe_n}}\; \ulk{\Dl'(\dDe_n)}\ru{6.5}$ by \bfcite{Ku05}{Prop.\ 2.6}. Therefore, the commutative quadrangle
\[
\xymatrix@C=13mm{
\ulk{\Dl(\dDe_n)}\ar[r]^{\ulk{G(\dDe_n)}}                                   & \ulk{\Dl'(\dDe_n)} \\
\ulk{\Dl^+(\b\De_n^\#)}\ar[u]^{(-)|_{\dDe_n}}\ar[r]^{\ulk{G^+(\b\De_n^\#)}} & \ulk{\Dl'^+(\b\De_n^\#)}\ar[u]_{(-)|_{\dDe_n}} \\
}
\]
shows that $\ulk{G^+(\b\De_n^\#)}$ is $1$-epimorphic; concerning notation, cf.\ \bfcite{Ku05}{\S1.2.1.1, \S1.2.1.3}. Therefore, we may define a transformation $\theta'_n$ for $\Dl'$ by the
requirement that 
\[
        \xymatrix{
        \ulk{\Dl^+(\b\De_n^\#)} \ar[rrr]^{\ulk{G^+(\b\De_n^\#)}}
              \ar@/_5mm/[dd]_{[-]^{+1}}="phf" 
              \ar@/^5mm/[dd]^{[-^{+1}]}="phg" 
        & & & \;\ulk{\Dl'^+(\b\De_n^\#)} 
              \ar@/_5mm/[dd]_{[-]^{+1}}="psf"
              \ar@/^5mm/[dd]^{[-^{+1}]}="psg" \\
        & & & & \\
        \ulk{\Dl^+(\b\De_n^\#)} \ar[rrr]^{\ulk{G^+(\b\De_n^\#)}} & & & \;\ulk{\Dl'^+(\b\De_n^{\,\#})}
        \ar@2 "phf"+<6.0mm,0mm>;"phg"+<-6.0mm,0mm>^{\theta_n}
        \ar@2 "psf"+<6.0mm,0mm>;"psg"+<-6.0mm,0mm>^{\theta'_n}
        }
\]
be commutative, i.e.\ that $\theta_n\st\ulk{G^+(\b\De_n^\#)} = \ulk{G^+(\b\De_n^\#)}\st\theta'_n\,$. In other words, there exists a unique $\theta'_n$ making this diagram commutative.

Let $\theta' := (\theta'_n)_{n\ge 0}\,$, where for $n = 0$, we make use of $\ulk{\Dl'^+(\b\De_0^{\,\#})} = 0$. We claim that $\theta'$ is a Heller triangulation on $(\Dl',\TTT')$, i.e.\ 
that $(\Dl',\TTT',\theta')$ is a Heller triangulated category. Once this is proven, we see that by construction, $\Dl\lraa{G}\Dl'\ru{5.0}$ is strictly exact; cf.\ \bfcite{Ku05}{Def.\ 1.5.(iii)}.

Suppose given $m,\, n\,\ge\, 0$ and a periodic monotone map $\b\De_n\llaa{p}\b\De_m\,$. To prove that \mb{$\ul{p}^\#\st\theta'_m \sollgl \theta'_n\st\ul{p}^\#$}, we may precompose with the $1$-epimorphic functor 
$\ulk{G^+(\b\De_n^\#)}$ to obtain
\[
\hspace*{-2mm}
\ba{rclclcl}
\ulk{G^+(\b\De_n^\#)}\st\ul{p}^\#\st\theta'_m \hspace*{-8mm}
& =                                                                          & \ul{p}^\#\st\ulk{G^+(\b\De_m^\#)}\st\theta'_m & = & \ul{p}^\#\st\theta_m\st\ulk{G^+(\b\De_m^\#)} 
                                                                                                                             &   &                                     \vspace*{1mm} \\
& \underset{\text{\scr Heller triangulated}}{\overset{(\Dl,\TTT,\theta)}{=}} & \theta_n\st\ul{p}^\#\st\ulk{G^+(\b\De_m^\#)}  & = & \theta_n\st\ulk{G^+(\b\De_n^\#)}\st\ul{p}^\#   
                                                                                                                             & = & \ulk{G^+(\b\De_n^\#)}\st\theta'_n\st\ul{p}^\# \; .\\
\ea
\]

Suppose given $n\,\ge\, 0$. To prove that \mb{$\ul{\ffk}_n\st\theta'_{n+1} \sollgl \theta'_{2n+1}\st\ul{\ffk}_n\,\,$}, we may precompose with the $1$\nbd-epimorphic functor 
$\ulk{G^+(\b\De_{2n+1}^\#)}$ to obtain
\[
\ba{rclcl}
\ulk{G^+(\b\De_{2n+1}^\#)}\st\ul{\ffk}_n\st\theta'_{n+1} \hspace*{-8mm}
& =                                                                          & \ul{\ffk}_n\st\ulk{G^+(\b\De_{n+1}^\#)}\st\theta'_{n+1}   &   &  \vspace*{1mm}\\
& =                                                                          & \ul{\ffk}_n\st\theta_{n+1}\st\ulk{G^+(\b\De_{n+1}^\#)}    &   &  \vspace*{1mm}\\
& \underset{\text{\scr Heller triangulated}}{\overset{(\Dl,\TTT,\theta)}{=}} & \theta_{2n+1}\st\ul{\ffk}_n\st\ulk{G^+(\b\De_{n+1}^\#)}   &   &  \vspace*{1mm}\\
& =                                                                          & \theta_{2n+1}\st\ulk{G^+(\b\De_{2n+1}^\#)}\st\ul{\ffk}_n  & = & \ulk{G^+(\b\De_{2n+1}^\#)}\st\theta'_{2n+1}\st\ul{\ffk}_n \; .\\
\ea
\]
\qed

\begin{Proposition}
\label{PropL5}
Recall that $(\Cl,\TTT,\tht)$ is a closed Heller triangulated category, and that $\Nl$ is a thick subcategory of $\Cl$.

There exists a unique Heller triangulation $\theta$ on $(\CNl,\TTT)$ such that $\Cl\lraa{\LLL}\CNl$ is strictly exact; cf.\ \mb{\bfcite{Ku05}{Def.\ 1.5}}.

Then $(\CNl,\TTT,\theta)$ is a closed Heller triangulated category; cf.\ Definition~{\rm\ref{DefClosed}.}
\end{Proposition}

{\it Proof.} By Lemma \ref{LemL2}, the category $\CNl$ is weakly abelian, and $\LLL : \Cl\lra\CNl$ is subexact. By Remark \ref{RemL3}, $\CNl$ 
carries a shift automorphism, and $\LLL$ is compatible with the shift automorphisms on $\Cl$ and on $\CNl$. By Lemma \ref{LemA1}, the functor 
$\Cl(\dDe_n)\lrafl{25}{\LLL(\dDe_n)} (\CNl)(\dDe_n)$ is $1$-epimorphic for $n\ge 0$.
Therefore, existence and uniqueness of $\theta$ follow by Lemma \ref{LemL4}.

It remains to be shown that $\CNl$ is closed. By isomorphic substitution, it suffices to show that each morphism in the image of $\LLL$ has a cone in $\CNl\,$; cf.\ \bfcite{Ku05}{Lem.~3.4.(6)}. But this follows
from $\Cl$ being closed and from $\LLL$ being strictly exact.
\qed

\bq
An object $(X\lraa{x} X')$ of the Freyd category $\h\Cl$ is called {\it $\Nl\!$-zero} if $x$ factors over an object of $\Nl\,$; concerning $\h\Cl$, cf.\ \bfcite{Ku05}{\S A.6.3}.
Note that an object of $\h\Cl$ that is isomorphic to a summand of an $\Nl\!$-zero object is itself $\Nl\!$-zero.

\begin{Remark}
\label{RemL6}
A morphism in $\,\Cl$ is an \Niso\ if and only if its kernel and its cokernel, taken in $\,\h\Cl$, are $\Nl\!$-zero.
\end{Remark}

Note that this criterion {\sf does not make reference to the Heller triangulated structure on $\Cl$,} but only to the fact that $\Cl$ is weakly abelian. One might ask for conditions on $\Nl$
that only use weak abelianess of $\Cl$, and that nonetheless suffice to turn $\Cl_{M(\Nl)}$ into a weakly abelian category -- where now $M(\Nl)$ is the subset of morphisms of $\Cl$
defined by the criterion given in Remark \ref{RemL6}.

{\it Proof of Remark~{\rm\ref{RemL6}.}} Suppose that $X\lraa{f} Y$ is an \Niso\ in $\Cl$. Then it has a weak kernel $N$ and a weak cokernel $M$ in $\Ob\Nl$. By construction of the kernel in $\h\Cl$, it is of the 
form $(N\lra X)$. Dually, the cokernel is of the form $(Y\lra M)$; cf.\ \bfcite{Ku05}{\S A.6.3}.

Conversely, suppose that the kernel and the cokernel of the morphism $X\lraa{f} Y$, taken in $\h\Cl$, are $\Nl\!$-zero. Consider the exact functor $\h\Cl\lraa{\h\LLL}(\Cl\fbby\Nl)\h\,$ that prolongs $\LLL$ on the
level of Freyd categories. It maps $f$ to an isomorphism, since in the abelian category $(\Cl\fbby\Nl)\h\,$, the image of $f$ has zero kernel and zero cokernel. Since 
$\Cl\fbby\Nl\lra (\Cl\fbby\Nl)\h\,$ is full and faithful, the image of $f$ under $\LLL$ in $\Cl\fbby\Nl$ is an isomorphism, too. Hence $f$ is an \Niso\ in $\Cl\,$; cf.\ Remark \ref{RemA0}.\qed
\eq

\begin{Proposition}[universal property]
\label{PropL7}
\hsp{3} Recall that $(\Cl,\TTT,\tht)$ is a closed Heller triangulated category, and that $\Nl$ is a thick subcategory of $\Cl$. 

Let $\theta$ be the unique Heller triangulation on $(\CNl,\TTT)$ such that the localisation functor $\Cl\lraa{\LLL}\CNl$ is strictly exact;
cf.\ Proposition~{\rm\ref{PropL5}.} Suppose given a Heller triangulated category $(\Cl',\TTT',\tht')$.

Recall that we write $\,\bo\,\Cl,\Cl'\bc_\text{\rm ex}$ for the category of exact functors and periodic transformations from $\,\Cl$ to $\,\Cl'\,$; cf.\ Definition~{\rm\ref{DefCatExFun}.}

Write $\,\bo\,\Cl,\Cl'\bc_\text{\rm ex, $\Nl$} \tm \bo\,\Cl,\Cl'\bc_\text{\rm ex}$ for the full subcategory consisting of exact functors $(F,a)$ such that $NF\iso 0$ for all $N\in\Ob\Nl$.

Recall that we write $\,\bo\,\Cl,\Cl'\bc_\text{\rm st\,ex}$ for the category of strictly exact functors and periodic transformations from $\,\Cl$ to $\,\Cl'\,$; cf.\  Definition~{\rm\ref{DefCatExFun}.}

Write $\,\bo\,\Cl,\Cl'\bc_\text{\rm st\,ex, $\Nl$} \tm \bo\,\Cl,\Cl'\bc_\text{\rm st\,ex}$ for the full subcategory consisting of strictly exact functors $F$ such that $NF\iso 0$ for all $N\in\Ob\Nl$.

\begin{itemize}
\item[\rm (1)] We have a strictly dense equivalence 
\[
\barcl
\bo\,\Cl,\Cl'\bc_\text{\rm ex, $\Nl$}            & \llafl{30}{\LLL\st\, (-)} & \bo\,\CNl,\,\Cl'\bc_\text{\rm ex} \\
(\LLL\st\,G\,,\,\LLL\st\,b)\= (\LLL,1)\st\,(G,b) & \llamapsfl{}{}            & (G,b) \; .                        \\
\ea
\]
\item[\rm (2)] We have a strictly dense equivalence 
\[
\barcl
\bo\,\Cl,\Cl'\bc_\text{\rm st\,ex, $\Nl$} & \llafl{30}{\LLL\st\, (-)} & \bo\,\CNl,\,\Cl'\bc_\text{\rm st\,ex} \\
\LLL\st\,G                                & \llamapsfl{}{}            &  G \; .                               \\
\ea
\]
\end{itemize}
\end{Proposition}

{\it Proof.} 

{\it Ad {\rm (1).}} Welldefinedness of the functor $\LLL\st\, (-)$ follows from $\LLL$ being strictly exact and exact functors being stable under composition; cf.\ Proposition~\ref{PropL5}, Remark~\ref{RemD1}.

We make use of the universal property of the localisation to the extent stated in Remark~\ref{RemA2}.

Suppose given exact functors $\Cl\,\lradoublea{(F,a)}{(G,b)}\,\Cl'\ru{-3}$ and a periodic transformation $F\lraa{u} G$.

Let $\br F : \CNl \lra \Cl'$ be defined by $\LLL\st\,\br F := F$. Let $\br G : \CNl \lra \Cl'$ be defined by $\LLL\st\,\br G := G$.

Recall that the shift on $\CNl$ is, abusively, also denoted by $\TTT$, so that $\TTT\st\LLL = \LLL\st\TTT$.
Let the transformations $\br a$ and $\br b$\ru{4.5} be defined by 
\[
\barcl
\LLL\st\,(\TTT\st\,\br F \lraa{\br a} \br F\st\TTT') & := & (\TTT\st\, F \lraa{a} F\st\TTT') \vsp{1} \\
\LLL\st\,(\TTT\st\,\br G \lraa{\br b} \br G\st\TTT') & := & (\TTT\st\, G \lraa{b} G\st\TTT')\; .     \\
\ea
\]
Let the transformation $\br{F}\lraa{\br u} \br{G}$ be defined by 
\[
\LLL\st\,(\br{F}\lraa{\br u} \br{G}) \;:=\; (F\lraa{u} G)\; . 
\]
We have to show that $(\br F,\br a)$ is exact and that $\br u$ is periodic.

{\it Ad $\br{F}$ exact.} Since $X\br{a} = X\LLL\br{a} = Xa$ is an isomorphism for $X\in\Ob\CNl = \Ob\Cl$, the transformation $a$ is an isotransformation.

To show that $\br{F}$ is subexact, by Lemma~\ref{LemC0_5}, it suffices to show that given a morphism $f$ in $\CNl$, it has a weak cokernel that is preserved by $\br{F}$.
By isomorphic substitution, we may assume that $f = f'\LLL$ for some morphism $f'$ in $\Cl$. Let $(f',g',h')$ be a $2$-triangle in $\Cl\,$; cf.\ Lemma~\ref{LemCl3}. 
Since $\LLL$ is strictly exact, the $2$-triangle $(f,g'\LLL,h'\LLL)$ results. In particular, $g'\LLL$ is a weak cokernel of $f$. Since $F$ is subexact, $g'F = g'\LLL\br F$ is a weak cokernel of 
$f'F = f'\LLL\br F = f\br F$. 

Suppose given $n\ge 0$. We shall make use of the abbreviation $\ulk{F} = \ulk{F^+(\bDes{n})}\,$, etc. It remains to show that
\[
(\theta_n\st \ulk{\br F})\cdot\ulk{\br a} \;\sollgl\; \ulk{\br F}\st\tht'_n\;.
\]
Since $\ul{\LLL} = \ulk{\LLL^+(\b\De_n^\#)}$ is $1$-epimorphic by Lemmata~\ref{LemA1}~and~\ref{LemL4}, it suffices to show that 
\[
\ul{\LLL}\st((\theta_n\st \ulk{\br F})\cdot\ulk{\br a}) \;\sollgl\; \ul{\LLL}\st\ulk{\br F}\st\tht'_n\;.
\]

In fact,
\[
\barcl
\ul{\LLL}\st ((\theta_n\st \ulk{\br F})\cdot\ulk{\br a})
& =                   & (\ul{\LLL}\st\theta_n\st \ulk{\br F})\cdot(\ul{\LLL}\st\ulk{\br a}) \vsp{1}\\
& \aufgl{$\LLL$ ex.}  & (\tht_n\st\ul{\LLL}\st\ulk{\br F})\cdot(\ul{\LLL}\st\ulk{\br a})    \vsp{1}\\
& =                   & (\tht_n\st\ulk{F})\cdot\ulk{a}                                      \vsp{1}\\
& \aufgl{$(F,a)$ ex.} & \ulk{F}\st\tht'_n                                                   \vsp{1}\\
& =                   & \ul{\LLL}\st\ulk{\br F}\st\tht'_n\; .                                      \\
\ea
\]

{\it Ad $\br u$ periodic.} We have to show that 
\[
(\TTT\st\,\br u)\cdot \br b \;\sollgl\; \br a\cdot (\br u\st\TTT')
\]
as transformations from $\TTT\st\,\br F$ to $\br G\st\TTT'$. By Remark~\ref{RemA2}, it suffices to show that
\[
\LLL\st\,((\TTT\st\,\br u)\cdot \br b) \;\sollgl\; \LLL\st\,(\br a\cdot (\br u\st\TTT'))\; . 
\]
In fact,
\[
\barcl
\LLL\st\,((\TTT\st\,\br u)\cdot \br b)
& =                & (\LLL\st\TTT\st\,\br u)\cdot (\LLL\st\,\br b)   \vsp{1}\\
& =                & (\TTT\st\LLL\st\,\br u)\cdot (\LLL\st\,\br b)   \vsp{1}\\
& =                & (\TTT\st\,u)\cdot b                             \vsp{1}\\
& \aufgl{$u$ per.} & a \cdot (u\st\TTT')                             \vsp{1}\\
& =                & (\LLL\st\,\br a) \cdot (\LLL\st\,\br u\st\TTT') \vsp{1}\\
& =                & \LLL\st\,(\br a \cdot (\br u\st\TTT')) \; .            \\
\ea
\]

{\it Ad {\rm (2).}} Welldefinedness of the functor $\LLL\st\,(-)$ follows from $\LLL$ being strictly exact and strictly exact functors being stable under composition; cf.\ Proposition~\ref{PropL5}, Remark~\ref{RemD1}.

Keep the notation of the proof of~(1).
Given an exact functor $(F,a)$ from $\Cl$ to $\Cl'$, we infer from $a = 1$, using $\LLL\st\,\br a = a$, that $\br a = 1$.
\qed


\appendix

\section{Some general assertions}
\label{SecGen}

\begin{footnotesize}

\bq
This appendix serves as a tool kit consisting of known results and folklore lemmata. We do not claim originality.
\eq

\subsection{Remarks on coretractions and retractions}
\label{SecAppRemAdd}

\begin{Remark}
\label{RemAdd1}
Let $\Al$ be a category. 

Suppose given $X$, $Z$ in $\Ob\Al$, and morphisms $X\lraa{i} Z\lraa{p} X$ such that $ip = 1_X\,$.

Suppose given $Y$, $W$ in $\Ob\Al$, and morphisms $Y\lraa{j} W\lraa{q} Y$ such that $jq = 1_Y\,$.

Suppose given $X\lraa{u} Y$ in $\Al$. Let $Z\lraa{v} W$ be defined by $v := puj$. Then $vq = pu$ and $iv = uj$.
\[
\xymatrix{
Z\ar[r]^p\ar[d]_v & X\ar[r]^i\ar[d]_u & Z\ar[d]^v \\
W\ar[r]^q         & Y\ar[r]^j         & W         \\
}
\]
\end{Remark}

{\it Proof.} We have $vq = pujq = pu$ and $iv = ipuj = uj$.\qed

\begin{Remark}
\label{RemAdd2}
Let $\Al$ be a category. Suppose given $Z,\,X,\,Z',\,W,\,Y,\,W'\,\in\,\Ob\Al$.

Suppose given morphisms $X\lraa{i} Z\lraa{p} X$ such that $ip = 1_X\,$.

Suppose given morphisms $X\lraa{i'} Z'\lraa{p'} X$ such that $i'p' = 1_X\,$.

Suppose given morphisms $Y\lraa{j} W\lraa{q} Y$ such that $jq = 1_Y\,$.

Suppose given morphisms $Y\lraa{j'} W'\lraa{q'} Y$ such that $j'q' = 1_Y\,$.

Suppose given $Z\lraa{v} W$ and $Z'\lraa{v'} W'$ such that $pi'v' = vqj'$. 

Then there exists a unique morphism $X\lraa{u} Y$ in $\Al$ such that $vq = pu$ and $i'v' = uj'$.
\[
\xymatrix{
Z\ar[r]^p\ar[d]_v & X\ar[r]^{i'}\ar[d]_u & Z'\ar[d]^{v'} \\
W\ar[r]^q         & Y\ar[r]^{j'}         & W'            \\
}
\]
If $v$ and $v'$ are isomorphisms, so is $u$.
\end{Remark}

{\it Proof.} Uniqueness follows from $p$ being epic and $j'$ being monic.

For existence, we let $u := ivq = i'v'q'$, the latter equality holding because of $pivqj' = pipi'v' = pi'v' = vqj' = vqj'q'j' = pi'v'q'j'$, using $p$ epic and $j'$ monic.
Then $pu = pi'v'q' = vqj'q' = vq$ and $uj' = ivqj' = ipi'v' = i'v'$.

If $v$ and $v'$ are isomorphisms, then let $u' := j v^- p = j' v'^- p'$ to get $uu' = ivqj' v'^- p' = ipi'v'v'^- p' = 1$ and $u'u = jv^-p\,i'v'q' = j v^- v q j' q' = 1$, so that $u' = u^-$.
In particular, $u$ is an isomorphism.
\qed

\subsection{Two lemmata on subexact functors}
\label{SecAppLemSubexact}

Suppose given weakly abelian categories $\Al$ and $\Al'\,$; cf.\ e.g.\ \bfcite{Ku05}{Def.~A.26.(3)}. Suppose given an additive functor $F : \Al\lra\Al'$. Recall that $F$ is called subexact
if the induced functor $\h F : \h\Al\lra\h\Al'$ on the Freyd categories is exact; cf.\ \bfcite{Ku05}{\S1.2.1.3}.

\begin{Lemma}
\label{LemC0_5}
The following assertions {\rm (1,\,2,\,3,\,3$^\0$,\,4,\,4$^\0$)} are equivalent.

\begin{itemize}
\item[\rm (1)] The functor $F$ is subexact.
\item[\rm (2)] The functor $F$ preserves weak kernels and weak cokernels.
\item[\rm (3)] The functor $F$ preserves weak kernels.
\item[\rm (3$^\0$)] The functor $F$ preserves weak cokernels.
\item[\rm (4)] For each morphism $X\lraa{t} Y$ in $\Al$, there exists a weak kernel $W\lraa{w} X$ such that $wF$ is a weak kernel of~$tF$.
\item[\rm (4$^\0$)] For each morphism $X\lraa{t} Y$ in $\Al$, there exists a weak cokernel $Y\lraa{w'} W'$ such that $w'F$ is a weak cokernel of~$tF$.
\end{itemize}
\end{Lemma}

{\it Proof.} Ad (1) $\Rightarrow$ (4). 
Suppose given a morphism $X\lraa{t} Y$ in $\Al$. Let $K\lramonoa{i} X$ be a kernel of~$t$ in $\h\Al$. 
Choose $A\lraepifl{30}{b} K$ with $A\in\Ob\Al$. Since $\h F$ is exact, $A\h F\lrafl{24}{(bi)\h F} X\h F \lrafl{24}{t\h F} Y\h F$
is exact at $X\h F$. So $(bi)\h F = (bi)F$ is a weak kernel of $t\h F = tF$ in $\Al'$.

Ad (4) $\Rightarrow$ (3). Given a morphism $X\lraa{t} Y$ in $\Al$ and a weak kernel $W\lraa{w} X$, a morphism $V\lraa{v} X$ is a weak kernel of $t$ if and only if both ($w$ factors over $v$) and ($v$ factors over $w$).
So if $wF$ is a weak kernel of $tF$, so is $vF$. Consequently, if $F$ preserves a single weak kernel of $t$, it preserves all of them. 

Ad (3) $\Rightarrow$ (2). This follows by \bfcite{Ku05}{Rem.~A.27}.

Ad (2) $\Rightarrow$ (1). Using duality and uniqueness of the kernel up to isomorphism, it suffices to show that $\h F$ maps a chosen kernel of a given morphism to a kernel of its image 
under $\h F$. Since $F$ preserves weak kernels,
this follows by construction of a kernel; cf.\ e.g.\ \mb{\bfcite{Ku05}{\S A.6.3, item (1) before Rem.~A.27}}.
\qed

\begin{Lemma}
\label{LemAdjSubex}
Suppose that $\Cl\lraa{F}\Cl'$ is subexact. Suppose given a functor $\Cl\llaa{G}\Cl'$.

\begin{itemize}
\item[\rm (1)] If $G\adj F$, then $G$ is subexact.
\item[\rm (1$^\0$)] If $F\adj G$, then $G$ is subexact.
\end{itemize}
\end{Lemma}

{\it Proof.} Ad (1). As an adjoint functor between additive categories, $G$ is additive.

Let $1\lraa{\eps} GF$ be a unit and $FG\lraa{\et} 1$ a counit of the adjunction $G\adj F$.

By Lemma~\ref{LemC0_5}, it suffices to show that $G$ preserves weak cokernels. Suppose given $X'\lraa{u} X\lraa{v} X''$ such that $v$ is a weak cokernel of $u$. We have to show that $Gv$ is a 
weak cokernel of $Gu$. Suppose given $t : XG\lra T$ such that $uG\cdot t = 0$. Then 
\[
u\cdot X\eps\cdot tF \= X'\eps\cdot uGF \cdot tF \= X'\eps\cdot (uG\cdot t)F \= 0\; .
\]
Since $v$ is a weak cokernel of $u$, we obtain a morphism $s : X'' \lra TF$ such that $v\cdot s = X\eps\cdot tF$. Then
\[
vG\cdot (sG\cdot T\et) \= X\eps G\cdot tFG\cdot T\et \= X\eps G\cdot XG\et\cdot t \=  t\; .
\]
\qed

\subsection{Karoubi hull}
\label{SecAppKarHull}

The construction of the Karoubi hull is due to {\sc Karoubi}; cf.\ \bfcite{Ka69}{III.II}.

Suppose given an additive category $\Al$. The {\it Karoubi hull} $\w\Al$ has 
\[
\Ob\w\Al \; := \; \{\, (A,e) \,:\, \text{$A\in\Ob\Al$, $e\in\liu{\Al}{(A,A)}$ with $e^2 = e$}\,\}
\]
and, given $(A,e),\,(B,f)\,\in\,\Ob\w\Al$, 
\[
\liu{\Al}{\big((A,e),\,(B,f)\big)} \; := \; \{\,u\in\liu{\Al}{(A,B)} \, : \, e\cdot u\cdot f = u\,\}\; .
\]
Then $\w\Al$ is an additive category, in which all idempotents are split.

Composition is inherited from $\Al$. We have a full and faithful additive functor
\[
\barcl
\Al           & \lraa{\KKK} & \w\Al                             \\
(X\lraa{u} Y) & \lramaps    & \big((X,1)\lraa{u} (Y,1)\big)\;,  \\
\ea
\]
which we often consider as an inclusion of a full subcategory.

Suppose given an additive category $\Bl$ in which all idempotents are split.

\begin{Remark}
\label{RemAppKar1}\rm
Write $\fbo\Al,\Bl\fbc_\text{add}$ for the category of additive functors and transformations between such from $\Al$ to~$\Bl$.
The induced functor $\fbo\Al,\Bl\fbc\llafl{27}{\KKK\st\,(-)} \fbo\w\Al,\Bl\fbc$ restricts to a strictly dense equivalence
\[
\fbo\Al,\Bl\fbc_\text{add}\;\llafl{27}{\KKK\st\,(-)}\; \fbo\w\Al,\Bl\fbc_\text{add}\; .
\]
\end{Remark}

\begin{Lemma}
\label{LemAppKar2}\rm
Suppose given an additive functor $\Al\lraa{I}\Al'$ to an additive category $\Al'$ in which all idempotents split. 
By Remark~\ref{RemAppKar1}, we obtain a functor $J : \w\Al\lra\Al'$, unique up to isomorphism, such that the following triangle of functors commutes.
\[
\xymatrix{
\Al\ar[d]_{\KKK}\ar[r]^I & \Al' \\
\w\Al\ar[ur]_J           &      \\
}
\] 
If $I$ is full and faithful, and if every object of $\Al'$ is a direct summand of an object in the image of $I$, then $J$ is an equivalence.
\end{Lemma}

By abuse of notation, in the situation of Lemma~\ref{LemAppKar2}, we also write $\w\Al = \Al'$ and consider $I$ to be an inclusion of a full subcategory.

\subsection{Multiplicative systems}
\label{SubsecMultSyst}

The construction of the quotient category of a Verdier triangulated category is due to {\sc Verdier}; cf.~\bfcit{Ve63}.

Suppose given a category $\Cl$. 

\begin{Definition}
\label{DefAppMult}\rm
A set $M$ of morphisms of $\Cl$ is called a {\it multiplicative system} in $\Cl$ if (Fr\,1-4) are satisfied. An element of $M$ is called an $M$-isomorphism and
denoted by $X\imp Y$.

\begin{itemize}
\item[(Fr\,1)] Each identity in $\Cl$ is an $M$-isomorphism.
\item[(Fr\,2)] Suppose given $X\lraa{f} Y\lraa{g} Z\lraa{h} W$ in $\Cl$ such that $fg$ and $gh$ are $M$-isomorphisms. 

              Then $f$, $g$, $h\,$ and $f\cdot g\cdot h\,$ are $M$-isomorphisms.
\item[(Fr\,3)] Suppose given $X\flradoublea{f}{g} Y\ru{-4}$ in $\Cl$. There exists an $M$-isomorphism $s$ such that $sf = sg$ if and only if there exists an $M$-isomorphism $t$ such that $ft = gt$.
\item[(Fr\,4)] Given 
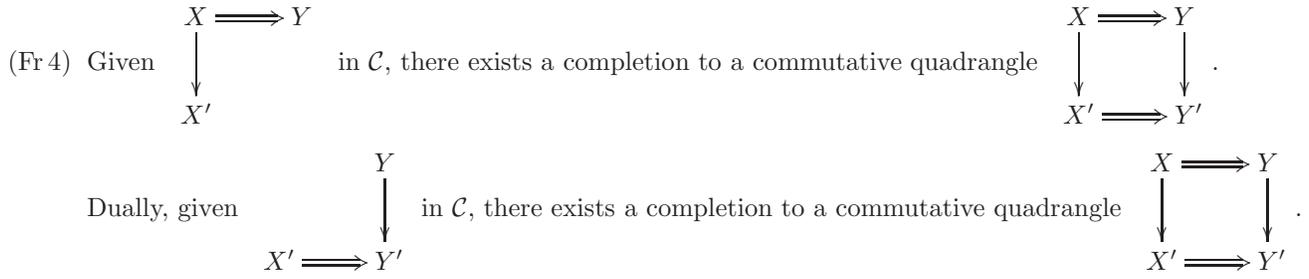
\begin{picture}(220,80)(0,-60)\put(0,0){
$
\xymatrix{
X\ar@{=>}[r]\ar[d]       & Y \\
X'                       &   \\
}
$
}\end{picture}\ru{-10}
in $\Cl$, there exists a completion to a commutative quadrangle
\begin{picture}(220,80)(0,-60)\put(0,0){
$
\xymatrix{
X\ar@{=>}[r]\ar[d]  & Y\ar[d] \\
X'\ar@{=>}[r]       & Y'\!\!  \\
}
$
}\end{picture}.

Dually, given
\begin{picture}(220,80)(0,-60)\put(0,0){
$
\xymatrix{
               & Y\ar[d] \\
X'\ar@{=>}[r]  & Y'\!\!  \\
}
$
}\end{picture}\ru{-10}
in $\Cl$, there exists a completion to a commutative quadrangle
\begin{picture}(220,80)(0,-60)\put(0,0){
$
\xymatrix{
X\ar@{=>}[r]\ar[d]        & Y\ar[d] \\
X'\ar@{=>}[r]             & Y'\!\!  \\
}
$
}\end{picture}.
\end{itemize}

Cf.\ \bfcite{Ve63}{\S 2, no.\ 1}.
\end{Definition}

Suppose given a multiplicative system $M$ in $\Cl$. Using (Fr\,2), we note that in the first assertion of (Fr\,4), if $X\lra X'$ is an $M$-isomorphism, then there exists a commutative completion with $Y\lra Y'$
being an $M$\nbd-isomorphism. And dually.

The category $\Cl_M\,$, called {\it localisation of $\Cl$ at $M$,} is defined as follows. Let $\Ob\Cl_M := \Ob\Cl$. A morphism from $X$ to $Y$ is a {\it double fraction,} which is 
an equivalence class of diagrams of the following form.
\[
\xymatrix@R=3mm{
  & X'\ar@{=>}[dl]_s\ar[r]^f & Y'\!\! &                 \\
X &                          &        & Y\ar@{=>}[ul]_t \\
}
\]

The diagrams $(s,f,t)$ and $(s's,s'ft',tt')$ are declared to be {\it elementarily equivalent}, provided $s'$ and $t'$
are $M$\nbd-isomorphisms. To form double fractions, we take the equivalence relation generated by elementary equi\-valence.

The equivalence class of the diagram $(s,f,t)$ is written $s\<\, f/t$. So $s\<\, f/t = \w s\<\, \w f/\w t$ if and only if there exist $M$-isomorphisms $u$, $\w u$, $v$, $\w v$ such that
$us = \w u\w s$ and $tv = \w t \w v$ and $ufv = \w u\w f \w v$. 
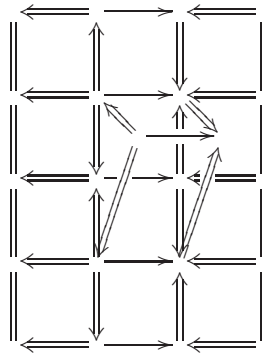
\[
\xymatrix@C-5mm@R-5mm{
\ar@{=}[dd] && \ar@{=>}[ll]\ar[rr]                         &                                   & \ar@{=>}[dd]                        && \ar@{=>}[ll]\ar@{=}[dd]             \\
            &&                                             &                                   &                                     &&                                     \\
\ar@{=}[dd] && \ar@{=>}[ll]\ar[rr]\ar@{=>}[uu]\ar@{=>}[dd] &                                   & \ar@{=>}[dr]                        && \ar@{=>}[ll]\ar@{=}[dd]             \\
            &&                                             & \ar@{=>}[ul]\ar@{=>}[dddl]\ar[rr] &                                     &&                                     \\
\ar@{=}[dd] && \ar@{=>}[ll]\ar[rr]|(0.33)\hole             &                                   & \ar@{=>}[uu]|(0.5)\hole\ar@{=>}[dd] && \ar@{=>}[ll]|(0.67)\hole\ar@{=}[dd] \\
            &&                                             &                                   &                                     &&                                     \\
\ar@{=}[dd] && \ar@{=>}[ll]\ar[rr]\ar@{=>}[uu]\ar@{=>}[dd] &                                   & \ar@{=>}[uuur]                      && \ar@{=>}[ll]\ar@{=}[dd]             \\
            &&                                             &                                   &                                     &&                                     \\
            && \ar@{=>}[ll]\ar[rr]                         &                                   & \ar@{=>}[uu]                        && \ar@{=>}[ll]                        \\
}
\]

Write $f/t := 1\<\, f/t$, called a {\it right fraction}, and $s\< f := s\<\, f/1$, called a {\it left fraction}. Using (Fr\,4), each morphism in $\Cl_M$ can be represented both by a left fraction and
by a right fraction. Given right fractions $f/t$ and
$\w f/\w t$, they are equal if there exist $M$-isomorphisms $u$, $v$ and $\w v$ such that $ufv = u\w f\w v$ and $tv = \w t\w v$. By (Fr\,3), this implies the existence of $M$-isomorphisms $v$, $\w v$ and 
$u'$ such that $f(vu') = \w f(\w vu')$ and $t(vu') = \w t(\w vu')$. Dually for left fractions. 

\sbq
 So double fractions are a self-dual way to represent morphisms in $\Cl_M\,$. Right or left fractions are more efficient in many arguments.
\seq

The composite of two double fractions $s\<\, f/t$ and $u\<\, g/v$ is defined, using (Fr\,4) for the commutative diagram
\[
\xymatrix@R=4mm{
  & \ar@{=>}[dl]_s\ar[r]^f &  \ar[r]^{g'} &                              &                         &                  &                       \\
X &                        &              & Y\!\!\ar@{=>}[ul]^t          &                         &                  & Z\ar@{=>}[dl]^v\zw{,} \\         
  &                        &              &  \ar[r]_{f'}\ar@{=>}[uull]^x & \ar@{=>}[ul]_u \ar[r]_g & \ar@{=>}[uull]_y &                       \\
}
\]
to be equivalently $s\<\, fg'/vy$ or $xs\<\, f'g/v$. By (Fr\,4,\,2,\,3), this definition is independent of the chosen completion with $g'$ and $y$, and, likewise, of the chosen completion with $x$ and $f'$.

Independence of the choice of the representative $s\<\, f/t$ is seen considering an elementary equivalence and using (Fr\,4,\,2), thus obtaining an elementary equivalence of the two possible representatives
of the composite. Likewise independence of the representative $u\<\, g/v$.

Associativity follows using right fractions and a commutative diagram constructed by means of (Fr\,4),
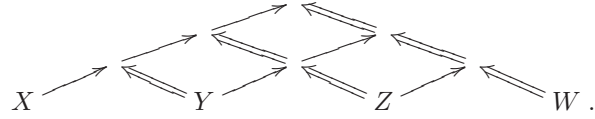
\[
\xymatrix@R=2mm{
         &         &                      &                      &                       &                   &                      \\
         &         & \ar[ur]              &                      & \ar@{=>}[ul]          &                   &                      \\
         & \ar[ur] &                      & \ar@{=>}[ul] \ar[ur] &                       & \ar@{=>}[ul]      &                      \\
X\ar[ur] &         & Y\ar@{=>}[ul]\ar[ur] &                      & Z \ar@{=>}[ul]\ar[ur] &                   & W\zw{.} \ar@{=>}[ul] \\
}
\]

Given $f\in\Mor\Cl$, we also write $1\<\, f/1 =: f$ in $\Cl_M$, by abuse of notation. Note that in $\Cl_M\,$, we have $s\<\, f/t = s^- f\, t^-$. 

\begin{Remark}
\label{RemA0}
A double fraction $s\<\, f/t$ represents an isomorphism in $\,\Cl_M$ if and only if $f$ is an $M$-isomorphism.
\end{Remark}

{\it Sketch.} First, using (Fr\,2), we reduce to the case of a right fraction $g/u$. For a right fraction in turn, the assertion follows applying (Fr\,2) to an associativity
diagram as above. \qed

\begin{Remark}
\label{RemA0_1}\rm
Given $X\flradoublea{f}{g} Y\ru{-4}$ in $\Cl$, we have $f = g$ in $\Cl_M$ if and only if there exists an $M$-isomorphism $t$ such that $ft = gt$ in $\Cl$, or, equivalently, if and only if there exists
an $M$-isomorphism $s$ such that $sf = sg$ in $\Cl$.
\end{Remark}

\begin{Remark}
\label{RemA0_3}
We have a functor $\Cl\lraa{\LLL}\Cl_M\,$, $f\lramaps 1\<\,f/1 = f$, called {\it localisation functor.} 

Given a category $\Tl$, we let $\fbo\,\Cl,\Tl\,\fbc_M$ be the full subcategory 
of $\fbo\,\Cl,\Tl\,\fbc$ consisting of functors that send all $M$\nbd-isomorphisms in $\Cl$ to isomorphisms in $\Tl$. The induced functor
\[
\fbo\,\Cl\,,\,\Tl\,\fbc_M \;\;\llafl{28}{\LLL\st\,(-)}\;\; \fbo\,\Cl_M\,,\,\Tl\,\fbc
\]
is a strictly dense equivalence, i.e.\ it is surjective on objects, full and faithful. 
\end{Remark}

{\it Sketch.} Given a functor $F\in\Ob\fbo\,\Cl\,,\,\Tl\,\fbc_M\,$, we may define $\br F$ on $\Cl_M$ by letting $X\br F := XF$ for $X\in\Ob\Cl_M = \Ob\Cl$ and by $(s\<\,f/t)\br F := (sF)^-\cdot (fF)\cdot (tF)^-$.
Then $\LLL\st\,\br F = F$. 

Given a transformation $(F\lraa{u} G)\in\Mor\fbo\,\Cl\,,\,\Tl\,\fbc_M\,$, we may define $\br F\lraa{\br u} \br G$ by setting $X\br u := Xu$ for \linebreak 
$X\in\Ob\Cl_M = \Ob\Cl$. Then $\LLL\st\,\br u = u$.
\qed

\begin{Lemma}
\label{LemA0_5}
Given $n\ge 0$, the functor
\[
\Cl(\dDe_n)\;\;\lrafl{25}{\LLL(\dDe_n)}\;\;\Cl_M(\dDe_n)\; ,
\]
given by pointwise application of $\,\LLL$, is dense.
\end{Lemma}

{\it Proof.} We may assume $n\ge 1$.

Suppose given $X\in\Ob\Cl_M(\dDe_n)$. To prove that for $i\in [1,n-1]$ there exists
an $X'\in\Ob\Cl_M(\dDe_n)$ isomorphic to $X$ such that $X'_j\lraa{x'} X'_{j+1}$ is in the image of $\LLL$ for $j\in [1,i-1]$, we proceed by induction on $i\ge 1$. Suppose the assertion to be true for $i$.
Let us prove the assertion for $i+1$.
Write $X'_{i}\lraa{x'} X'_{i+1}$ as a right fraction $f/s$. If $i = n-1$, we replace $X'_{i+1}$ by the target of $f$, and $f/s$ by $f$. If $i \le n-2$, we write $X'_{i+1}\lraa{x'} X'_{i+2}$ 
as a right fraction $g/u$ and construct the following commutative diagram using (Fr\,4).
\[
\xymatrix@!@R=-4mm{
                 &                       &                                 &                   &                        \\
                 & X''_{i+1}\ar[ur]^{g'} &                                 & \ar@{=>}[ul]_{s'} &                        \\
X'_{i}\ar[ur]^f  &                       & X'_{i+1}\ar@{=>}[ul]_s\ar[ur]^g &                   & X'_{i+2}\ar@{=>}[ul]_u \\
}
\]
We replacing the object $X'_{i+1}\ru{-2.5}$ by $X''_{i+1}\,$, the morphism $f/s$ by $f$ and the morphism $g/u$ by $g'/us'$. 

In both cases, we obtain a diagram isomorphic to $X'$ that conincides with $X'$ on $[1,i]$ and whose morphism from $i$ to $i+1$ is in the image of $\LLL$. \qed

\begin{Lemma}
\label{LemA1}
Given $n\ge 0$, the functor
\[
\Cl(\dDe_n)\;\;\lrafl{25}{\LLL(\dDe_n)}\;\;\Cl_M(\dDe_n)
\]
is $1$-epimorphic.
\end{Lemma}

{\it Proof.} We shall apply \bfcite{Ku05}{Lem.\ A.35}. By Lemma \ref{LemA0_5}, $\LLL(\dDe_n)$ is dense.

Suppose given $X,\, Y\,\in\,\Ob\Cl(\dDe_n)$ and a morphism $X\LLL(\dDe_n)\lraa{g} Y\LLL(\dDe_n)$ in $\Cl_M(\dDe_n)$. Let $g_i$ be represented by a right fraction $f_i/s_i$ for $i\in [1,n]$.

\hypertarget{Claim.App}{}%
We {\it claim} that for $i\in [1,n]$, we can find representatives $f'_j/s'_j$ for $j\in [1,i]$ such that there exist $h_j$ with $s'_j h_j = y s'_{j+1}$ and $f'_j h_j = x f'_{j+1}$ in $\Cl$ 
for $j\in [1,i-1]$. 
Let $f'_1 := f_1$ and $s'_1 := s_1\,$. Proceeding by induction on $i$, we have to write the right fraction $f_{i+1}/s_{i+1}$ suitably as $f'_{i+1}/s'_{i+1}\,$. 
First of all, by (Fr\,4), we find an $M$-isomorphism $\sa$ and a morphism $\xi$ such that $y\sa = s'_i\xi$ in $\Cl$. We have 
$f'_i\xi\sa^- = x f_{i+1} s_{i+1}^-$ in $\Cl_M\,$. Using (Fr\,4) and (Fr\,2), we find $M$-isomorphisms $s'$ and $\sa'$ such that $\sa s' = s_{i+1}\sa'$ in $\Cl$. Hence 
\[
f'_i\xi s' \= f'_i\xi s' s'^- \sa^- s_{i+1} \sa' \= f'_i\xi \sa^- s_{i+1} \sa' \= x f_{i+1} s_{i+1}^- s_{i+1} \sa' \= x f_{i+1} \sa'
\]
in $\Cl_M\,$. Composing with a further $M$-isomorphism, we may assume that $f'_i\xi s' = x f_{i+1} \sa'$ in $\Cl\,$; cf.\ Remark~\ref{RemA0_1}. 
We take $h_i := \xi s'$ and $s'_{i+1} := s_{i+1}\sa'$ and $f'_{i+1} := f_{i+1}\sa'$.
\[
\xymatrix{
X_i\ar[rrr]^x\ar[dd]_{f'_i}       &                   & & X_{i+1}\ar[dd]^{f_{i+1}}                        \\
                                  &                   & &                                                 \\
\ar[r]^\xi                        & \ar@{=>}[ur]^{s'} & & \ar@{=>}[ul]_{\sa'}                             \\ 
Y_i\ar@{=>}[u]^{s'_i}\ar[rrr]^y   &                   & & Y_{i+1}\ar@{=>}[ull]^{\sa}\ar@{=>}[u]^{s_{i+1}} \\ 
}
\]
This proves the \hyperlink{Claim.App}{\it claim}{\it,} in particular for $i = n$.
\[
\xymatrix{
X_1\ar[r]^x\ar[d]_{f'_1}      & X_2\ar[r]\ar[d]_{f'_2}      & \cdots\ar[r] & X_{n-1}\ar[r]^x\ar[d]_{f'_{n-1}}       & X_n\ar[d]^{f'_n}      \\
\ar[r]^{h_1}                  & \ar[r]                      & \cdots\ar[r] & \ar[r]^{h_{n-1}}                       &                       \\
Y_1\ar[r]^y\ar@{=>}[u]^{s'_1} & Y_2\ar[r]\ar@{=>}[u]^{s'_2} & \cdots\ar[r] & Y_{n-1} \ar[r]^y\ar@{=>}[u]^{s'_{n-1}} & Y_n\ar@{=>}[u]^{s'_n} \\
}
\]
Condition (C) of loc.\ cit.\ is satisfied letting the epizigzag have length $0$, letting the monozigzag be the single backwards diagram morphism consisting of the morphisms $s'_i\,$, and letting the required
diagram morphism in the image of $L(\dDe_n)$ consist of the morphisms $f'_i\,$. 
\qed

\begin{Remark}
\label{RemA2}
Suppose the category $\Cl$ to be additive. 

\begin{itemize}
\item[\rm (1)] An object $X$ is isomorphic to $0$ in $\,\Cl_M$ if and only if $\;X\imp 0$, or, equivalently, if and only if $\;0\imp X$. 
\item[\rm (2)] The category $\Cl_M$ is additive, and the functor $\LLL : \Cl\lra\Cl_M$ is additive. 
\item[\rm (3)] Given an additive category $\Tl$, the strictly dense equivalence
\[
\fbo\,\Cl\,,\,\Tl\,\fbc_M \;\;\llafl{28}{\LLL\st\,(-)}\;\; \fbo\,\Cl_M\,,\,\Tl\,\fbc
\]
restricts to a strictly dense equivalence from the category of additive functors from $\Cl_M$ to $\Tl$ to the category of additive functors from $\Cl$ to $\Tl$ that sends all $M$-isomorphisms to 
isomorphisms, written
\[
\fbo\,\Cl\,,\,\Tl\,\fbc_\text{\rm add,\,$M$} \;\;\llafl{28}{\LLL\st\,(-)}\;\; \fbo\,\Cl_M\,,\,\Tl\,\fbc_\text{\rm add}\; .
\]
\end{itemize}
\end{Remark}

{\it Sketch.} 

Ad~(1). If $X$ is isomorphic to $0$, then $X\imp X'\isimp 0\,$; cf.\ Remark~\ref{RemA0}. By (Fr\,4), we conclude that $0\imp X$.

Ad~(2). Given $X,\,Y\,\in\Ob\Cl$, the direct sum $X\ds Y$, together with $X\lrafl{28}{\smatez{1}{0}} X\ds Y$ and $Y\lrafl{28}{\smatez{0}{1}} X\ds Y$, remains a coproduct in $\Cl_M\,$. 

For existence of an induced morphism from the coproduct, we use (Fr\,4,\,2) to produce a common denominator of two right fractions. 

To prove uniqueness of the induced morphism, we suppose given 
$\smatze{f}{g}\!/s$ and $\smatze{f'}{g'}\!/s$, without loss of generality with common denominator, such that $f/s = \smatez{1}{0}\cdot \smatze{f}{g}\!/s = \smatez{1}{0}\cdot \smatze{f'}{g'}\!/s = f'/s$
and $g/s = \smatez{0}{1}\cdot \smatze{f}{g}\!/s = \smatez{0}{1}\cdot \smash{\smatze{f'}{g'}\!/s}\ru{-2} = g'/s$ in $\Cl_M\,$. So there exists an $M$-isomorphism $u$ such that $fu = f'u$, and an $M$-isomorphism 
$v$ such that $gv = g'v$, both in $\Cl$. By (Fr\,4,\,2), we obtain a common $M$-isomorphism $w$ such that $fw = f'w$ and $gw = g'w$ in $\Cl$. Hence $\smatze{f}{g} w = \smatze{f'}{g'} w$ in $\Cl$. Therefore
$\smatze{f}{g}\!/s = \smatze{f'}{g'}\!/s$ in $\Cl_M\,$.

Moreover, the automorphism $\smatzz{1}{0}{1}{1}$ of $X\ds X$ remains an automorphism in $\Cl_M\,$.

Ad~(3). Since $\LLL$ is additive, $\LLL\st\,(-)$ sends additive functors to additive functors. Conversely, given an additive functor $F : \Cl\lra\Tl$, the functor $\br F$ as constructed in the 
proof of Remark~\ref{RemA0_3} is additive.
\qed

\end{footnotesize}


\begin{footnotesize}

\parskip1.2ex

\vspace*{5mm}

\begin{flushright}
Matthias K\"unzer\\
University of Stuttgart\\
Fachbereich Mathematik\\
Pfaffenwaldring 57\\
D-70569 Stuttgart\\
kuenzer@mathematik.uni-stuttgart.de \\
w5.mathematik.uni-stuttgart.de/fachbereich/Kuenzer/\\
\end{flushright}
\end{footnotesize}


\end{document}